\documentclass[11pt]{article}
\usepackage{amsmath,amsfonts,amssymb,amsthm,graphicx}
\usepackage{dsfont}
\DeclareMathOperator{\sign}{sgn}
\usepackage{bbm}
\usepackage{enumerate}
\usepackage{stmaryrd}  
\usepackage{color}
\usepackage{natbib}
\usepackage[colorlinks=true,citecolor=blue,pdfpagemode=UseNone,pdfstartview=FitH]{hyperref}
\usepackage{tasks}
\usepackage{indentfirst}
\usepackage{multirow}

\newif\ifFULL
\ifFULL

\fi

\renewcommand{\d}{\,\mathrm{d}}

\newcommand{\p}{\mathbb P}

\newcommand{\laweq}{\overset{\rm d}{=}}

\newcommand{\R}{\mathbb{R}}
\newcommand{\N}{\mathbb{N}}
\newcommand{\id}{\mathds{1}}
\renewcommand{\ge}{\geqslant}
\renewcommand{\le}{\leqslant}
\renewcommand{\epsilon}{\varepsilon}

\theoremstyle{plain}
\newtheorem{theorem}{Theorem}

\newtheorem{proposition}{Proposition}
\theoremstyle{definition}
\newtheorem{definition}{Definition}
\newtheorem{example}{Example}

\theoremstyle{remark}
\newtheorem{remark}{Remark}
\theoremstyle{definition}

\usepackage[misc]{ifsym}

\begin{document}

\title{Trade-off between validity and efficiency of merging p-values under arbitrary dependence}
\author{Yuyu Chen\thanks{Department of Statistics and Actuarial Science,
  University of Waterloo.
  E-mail: \texttt{y937chen@uwaterloo.ca}.} \and Peng Liu\thanks{Department of Mathematical Sciences,
  University of Essex.
  E-mail: \texttt{peng.liu@essex.ac.uk}.}  \and Ken Seng Tan\thanks{Nanyang Business School, Nanyang Technological University.  E-mail: \texttt{kenseng.tan@ntu.edu.sg}.} \and
  Ruodu Wang\thanks{Department of Statistics and Actuarial Science,
  University of Waterloo.
  E-mail: \texttt{wang@uwaterloo.ca}.}
  }
 \date{}
\maketitle
\begin{abstract}
 Various methods of combining individual p-values into one p-value are widely used in many areas of statistical applications.
We say that a combining method is  valid for arbitrary dependence (VAD)  if it does not require any assumption on the dependence structure of the  p-values, whereas it is valid for some dependence (VSD) if it requires some specific, perhaps realistic but unjustifiable, dependence structures.
The trade-off between validity and efficiency of these methods is studied via analyzing  the choices of critical values under different dependence assumptions.
We introduce the notions of independence-comonotonicity balance (IC-balance)
and the price for validity.
In particular,
IC-balanced methods always produce an identical critical value for independent and perfectly positively dependent p-values,   a specific type of insensitivity to a family of  dependence assumptions.
We show that, among two very general classes of merging methods commonly used in practice, the  Cauchy combination method and the Simes method are the only IC-balanced ones.
Simulation studies and a real data analysis are conducted to analyze the sizes and powers of various combining methods in the presence of weak and strong dependence.
\end{abstract}

\newpage
\section{Introduction}

In many areas of statistical applications
where multiple hypothesis testing is involved,
the task of
merging several p-values into one naturally arises.
Depending on the specific application,
 these p-values may be from a single hypothesis or   multiple hypotheses, in small or large numbers, independent or correlated,  and with sparse or dense signals, leading to  different considerations when choosing merging procedures.

Let $K$ be a positive integer, and $F:[0,1]^{K}\rightarrow [0,\infty)$ be an increasing Borel function used to combine  $K$ p-values, which we shall refer to as a \emph{combining function}.
 Generally, the combined value may not be a valid p-value itself, and a critical point needs to be specified.  Different dependence assumptions on the p-values   lead  to significantly different critical points, and thus different
  statistical decisions.
The problem of merging p-values  has a long history, and early results can be found in \cite{T31}, \cite{pearson1933method} and \cite{fs1938methods} where   p-values are assumed to be independent.  Based on an idea of Tukey, \cite{donoho2004higher} developed the higher criticism statistics to detect weak and sparse signals effectively using independent p-values. Certainly, these methods do not always produce a valid p-value if  the assumption of independence is violated.
On the other hand, the independence assumption is often very difficult or impossible to verify in many applications where only one set of p-values is available.

There are, however, some methods that produce valid p-values without any dependence assumption.
A classic one is the Bonferroni method by taking the minimum of the p-values times $K$ (we allow combined p-values to be greater than $1$ and they can be treated as $1$) or equivalently, dividing the critical value by $K$.
 Other methods that are valid without assumptions
 include the ones based on order statistics  by \cite{ruger1978maximale} and \cite{hommel1983tests}, and the ones based on averaging by \cite{VW19};
   details of these merging methods are presented in Section \ref{s3}.

 Some other methods work under weak or moderate dependence assumptions, such as the   method of \cite{S86},
 which uses the minimum of  $Kp_{(i)}/i$ over $i=1,\dots, K$,
 where $p_{(i)}$ is the $i$-th smallest order statistic of $p_1,\dots,p_K$. 
The validity of the Simes method is shown under a large class of dependence structures  (e.g., \cite{sarkar1998some, sarkar2008simes}; \cite{benjamini2001control} and  \cite{rodland2006simes}), although even   such dependence assumptions are unlikely to hold in practice (see e.g., \citet[p.51]{E10}).
Two more recent methods include
the Cauchy combination test
proposed by \cite{LX19} using the weighted average of Cauchy transformed p-values,
and
the harmonic mean p-value of \cite{W19} using the harmonic mean of p-values.
Under mild dependence assumptions, these two methods are    asymptotically valid  as the significance level goes to $0$ (see Theorem \ref{theorem2}).

This paper is dedicated to a comprehensive and unifying treatment of p-value merging methods  under various dependence assumptions.
Some methods are
valid without any assumption  on the interdependence of p-values,
and they will be referred to as \emph{VAD methods}. On the other hand, methods that are valid for  some specific but realistic dependence assumption (e.g., independence, positive dependence, or joint normality dependence)
 will be referred to as \emph{VSD methods}.
Our main goal is to understand the difference and the trade-off between these methods.

For a fixed combining function $F$,
using a VAD method means choosing a smaller critical value (threshold) for making rejections compared to a VSD method. Thus, the gain of validity comes at the price of a loss of detection power.
As it is often difficult to make valid statistical inference on the dependence structure of p-values, our analysis also helps to understand the relative performance of VSD combining methods under the presence of model misspecification.
As a byproduct, we obtain several new theoretical results on the popular  Simes, harmonic, and Cauchy merging methods.

In the next section, we collect some basic definitions of VAD and VSD merging methods and their corresponding threshold functions.
We focus on symmetric merging functions for the tractability in their comparison. In Section \ref{s3}, we introduce two general  classes of combining functions, which include all methods mentioned above. Formulas for their VAD and VSD threshold functions are derived, some based on results from robust risk aggregation, e.g., \cite{WPY13}.
 In Section \ref{s4}, we introduce independence-comonotonicity balanced  (IC-balanced) combining functions, which are indifferent between the two dependence assumptions. We show that the Cauchy combination method and the Simes method are   the only IC-balanced ones among two general classes of combining methods, thus highlighting their unique roles. In Section \ref{sec:4p}, we establish strong similarity between the Cauchy combination and the harmonic averaging methods, and obtain an algebraic relationship between    the harmonic averaging and the Simes functions. In Section \ref{s5}, the {price for validity} is introduced to assess the loss of power of VAD methods compared to their VSD versions. Simulation studies and a real data analysis are conducted to analyze the relative performance of these methods. Simulation studies and a real data analysis are presented in Section \ref{sec:6} to analyze the relative performance of these methods. Proofs of all technical results are put in the supplementary material.

We conclude the section by providing additional notation and terminology that will be adopted in this paper.
All random variables are defined on an atomless probability space $(\Omega,\mathcal F,\p)$.  Random variables $X_1,\dots,X_n$ are comonotonic if there exist increasing functions $f_1,\dots,f_n$ and a random variable $Z$ such that $X_i=f_i(Z)$ for each $i=1,\dots,n$.
For $\alpha\in(0,1]$, $q_\alpha(X)$ is the left $\alpha$-quantile of a random variable $X$, defined as
\begin{align*}
    q_{\alpha}(X)&=\inf\{x\in \mathbb{R}\mid \mathbb{P}(X\leq x)\ge \alpha\}.
\end{align*}
We also use $F^{-1}(\alpha)$ for $q_\alpha(X)$ if $X$ follows the distribution $F$.
The set
 $\mathcal{U}$ is the set of  all  standard uniform random variables  defined on $(\Omega,\mathcal F,\p)$ (i.e., the set of all measurable functions on $(\Omega,\mathcal{F})$ whose distribution  under $\p$ is uniform on $[0,1]$)  and $\id$ is the indicator function.
 The equality $\laweq$ represents equality in distribution.
  For given $p_{1},\dots,p_{K}$, the order statistics  $p_{(1)},\dots,p_{(K)}$ are ordered from the smallest to the largest. The  equivalence $A_x \sim B_x$ as $x\to x_0$ means that $A_x/B_x\to 1$ as $x\to x_0$.
 All terms of ``increasing" and ``decreasing" are in the non-strict sense.

 \section{Merging methods and thresholds}\label{S2}
 Following the terminology of \cite{VW19},
a \emph{p-variable} is a random variable $P$ such that $\mathbb{P} (P\leq \epsilon) \leq \epsilon$, for all $\epsilon\in (0,1)$ (such random variables are called \emph{superuniform} by \cite{RBWJ19}). Values realized by p-variables are p-values. In the Introduction, p-values are used loosely for p-variables, which should be clear from the context.

Let $P_{1},\dots,P_{K}$ be $K$ p-variables for testing a common hypothesis.
A \emph{combining} \emph{function} is an increasing Borel measurable function $F:[0,1]^{K}\rightarrow [0,\infty)$ which transforms $P_{1},\dots,P_{K}$ into a single random variable $F(P_{1},\dots,P_{K})$. The choice of combining function depends on how one integrates   information, and some  common options are mentioned in the Introduction. Generally, $F(P_{1},\dots,P_{K})$ may not be a  valid p-variable.
For different choices of $F$ and assumptions on $P_1,\dots,P_K$,
one needs to assign a critical value $g(\epsilon)$ so that the hypothesis can be rejected with significance level $\epsilon\in(0,1)$ if $F(P_1,\dots,P_K)< g(\epsilon)$.
We call $g$ a \emph{threshold (function)} for $F$ and $P_1,\dots,P_K$.
Clearly, $g(\epsilon)$ is increasing in $\epsilon$.
In case $g$ is strictly increasing, which is the most common situation,
the above specification of $g$ is equivalent to requiring $g^{-1}\circ F(P_{1},\dots,P_{K})$ to be a p-variable. To objectively compare various combining methods,
one should compare the corresponding values of the function $g^{-1}\circ F$.

In some situations,
it might be convenient and  practical to assume
additional information on dependence structure of p-variables, e.g., independence, comonotonicity (i.e., perfectly positive dependence), and specific copulas. The choice of the threshold $g$ certainly depends on such assumptions. If no assumption is made on the interdependence of the p-variables, the corresponding threshold function is called a \emph{VAD threshold}, otherwise it is a \emph{VSD} \emph{threshold}.
 A testing procedure based on  a VAD threshold always produces a size  less than or equal to the significance level regardless of the dependence structure of the p-variables.

We denote the VAD threshold of a combining function $F$ by $a_F$.
If a merging method is valid for independent (resp.~comonotonic) dependence of p-variables, we use $b_{F}$ (resp.~$c_{F}$) to denote the corresponding valid threshold function, and we call it  the \emph{VI} (resp.~\emph{VC}) \emph{threshold}.
More precisely, for the equation
\begin{align}\label{e1}
 \mathbb{P}( F(P_{1},\dots,P_{K})<g(\epsilon)) \leq \epsilon,~~~\epsilon\in (0,1),
\end{align}
a VAD {threshold}  $g=a_F$ satisfies \eqref{e1}
for all p-variables $P_1,\dots,P_K$;
a VI {threshold}  $g=b_F$ satisfies \eqref{e1}
for all independent  p-variables $P_1,\dots,P_K$, and
 a VC {threshold}  $g=c_F$ satisfies \eqref{e1}
 for all comonotonic p-variables $P_1,\dots,P_K$.

 The comonotonicity assumption on the p-variables to combine (actually they are identical if they are uniform on $[0,1]$) is not interesting by itself for statistical practice. Nevertheless, comonotonicity is a benchmark for (extreme) positive dependence, and we analyze $c_F$ for the purpose of comparison; it helps us to understand how valid thresholds for different methods vary as the dependence assumption gradually shifts from independence to  extreme positive dependence. This  point will be made more clear in Sections \ref{s4}-\ref{s5}.

An immediate observation is that the p-variables  can be equivalently replaced by uniform random variables on $[0,1]$ as  for each p-variable $P$, we can find $U\in\mathcal U$ with $U\le P$;  see e.g., \cite{VW19}. 
Therefore, it suffices to consider p-variables  in $\mathcal U$.
Moreover, if $g$ satisfies \eqref{e1},
then any function that is smaller than $g$
is also valid. Hence, for the sake of power, it is natural to  use the largest functions that satisfy \eqref{e1}.
Putting these considerations together, we formally define the thresholds of interest  as follows.

\begin{definition}\label{def:1}
  The  thresholds $a_F$, $b_F$ and $c_{F}$ of a combining function $F$ are given by, for   $\epsilon\in (0,1)$,
  \begin{gather}
  a_{F}(\epsilon) =\inf\{q_{\epsilon}(F(U_{1},\dots,U_{K}) )\mid U_{1},\dots,U_{K}\in\mathcal{U}\},\label{e1p}\\
  b_{F}(\epsilon)=q_{\epsilon}(F(V_{1},\dots,V_{K}) ),\label{e2p}\\
  c_{F}(\epsilon)=q_{\epsilon}(F(U,\dots,U) ),\label{e3p}
  \end{gather}
where
  $U,V_{1},\dots,V_{K}$ are independent standard uniform random variables.

\end{definition}
{In what follows, we focus on the thresholds in Definition \ref{def:1}.} It is clear that $g=a_F$, $b_F$ or $c_F$  in Definition \ref{def:1} satisfies
\eqref{e1} under the respective dependence assumptions.

\begin{remark}
While the objects $b_F$ and $c_F$ in \eqref{e2p}-\eqref{e3p}  can often be explicitly calculated,
the object $a_F$ in \eqref{e1p} is generally  difficult to calculate for a chosen function $F$ due to the infimum taken over all possible dependence structures.
Techniques in the field of robust risk aggregation, in particular, results in \cite{WPY13}, \cite{embrechts2013model,EWW15} and \cite{W17},
are designed for such calculation, as illustrated by \cite{VW19}. By definition, for any threshold $g(\epsilon) >a_F(\epsilon)$, there exists some dependence structure of $(P_1,\dots,P_K)$ such that validity is lost, i.e., \eqref{e1} is violated. Moreover, if the combining function $F$ is continuous, the infimum in \eqref{e1p} is attainable; the proof of  this statement is similar to that of Lemma 4.2 of \cite{bernard2014risk}.
\end{remark}

\section{Combining functions}\label{s3}

\subsection{Two general classes of combining functions}
We first  introduce two general classes of combining functions, the generalized mean class and the order statistics class. 
 Let $p_{1},\dots,p_{K}\in [0,1]$ be the $K$ realized p-values. The first class of combining functions is the generalized mean, that is,
\begin{equation*}
    M_{\phi,K}(p_{1},\dots,p_{K})=\phi^{-1}\left(\frac{1}{K}\sum_{i=1}^{K}\phi(p_{i})\right),
\end{equation*}
where $\phi:[0,1]\rightarrow[-\infty,\infty]$ is a  continuous and strictly monotone function and $\phi^{-1}$ is its inverse on the domain $\phi([0,1])$. Many combining functions used in the statistical literature are included in this class. For example,
the Fisher method (\cite{fs1938methods}) corresponds to
the geometric mean with $\phi(p)=\log(p)$;
the averaging methods of \cite{VW19} and \cite{W19} correspond  to the functions $\phi(p)=p^{r}$, and $r\in[-\infty,\infty]$ (including  limit cases), and the Cauchy combination method of \cite{LX19}   corresponds to $\phi(p)=\tan\left(\pi\left(p-\frac{1}{2}\right)\right)$.

The second class of combining functions is built on order statistics.
Let $\alpha=\left(\alpha_{1},\dots,\alpha_{K}\right)\in\mathbb{R}^{K}_{+}$, where $\R_+=[0,\infty)$.  We define the combining function
\begin{equation*}
    S_{\alpha,K}(p_{1},\dots,p_{K})=\min_{i\in\{1,\dots,K\}}\frac{p_{(i)}}{\alpha_{i}},
\end{equation*}
 where the convention is $p_{(i)}/\alpha=\infty$ if $\alpha=0$.
 If $\alpha_{1}=1/K$ and all the other components of $\alpha$ are 0, then using $S_{\alpha,K}$ yields the Bonferroni method based on the minimum of p-values.
 The VAD method via order statistics of \cite{ruger1978maximale} uses $S_{\alpha,K}$ by setting $\alpha_{i}=i/K$ for a fixed $i\in\{1,\dots,K\}$ and  all the other components of $\alpha$ to be 0.
On the other hand,
if $\alpha_{i}=i/K$ for each $i =1,\dots,K$,
then  we arrive at the method of  \cite{S86}; in this case, we will simply denote $S_{\alpha,K}$ by $S_K$, namely,  $$S_K(p_1,\dots,p_K):=   \min_{i\in\{1,\dots,K\}}\frac{K p_{(i)}}{i}, $$ and  $S_K$ will be called  the \emph{Simes function}.
The method of  \cite{hommel1983tests} uses $ \ell_K S_K$, which is $S_K$ adjusted via the VAD threshold,
where
\begin{align}
\label{eq:ellK}
\ell_{K}=\sum_{k=1}^K\frac{1}{k}.
\end{align}
If $\alpha_{i+1} \le \alpha_i$,
then the term $p_{(i+1)}/\alpha_{i+1}$ does not contribute to the calculation of $ S_{\alpha,K}(p_{1},\dots,p_{K})$.
Hence, we can safely replace $\alpha_{i+1}$ by $\alpha_i$ without changing the function $ S_{\alpha,K}$.
Thus, we shall assume, without loss of generality, that $\alpha_1\le \dots\le \alpha_K$.
Admissibility of VAD merging methods in the above two classes are studied by \cite{VWW20}.

Recall that a function $F:\R_+^K\to \R$ is homogeneous if
$F(\lambda \mathbf x) =\lambda F(\mathbf x)$ for all $\lambda >0$ and $\mathbf x\in \R_+^K$.
It is clear that the function $S_{\alpha,K}$ is homogeneous, and so are  the averaging methods of \cite{VW19}. In such cases, we can show that  the VAD threshold $a_F$ is a linear function.

\begin{proposition}\label{p1}
If the combination function $F$ is homogeneous, then the VAD threshold $a_F(x)$ is a constant times $x$ on $(0,1)$.
\end{proposition}

In the subsections below we will discuss several special cases of the above two classes of combining functions, and analyze their corresponding threshold functions.
As the first example, we note that the functions $a_F$, $b_F$ and $c_F$ for the Bonferroni method can be easily verified.
\begin{proposition}\label{prop:bonf}
Let $F(p_1,\dots,p_K)=\min\{p_1,\dots,p_K\}$ for $p_1,\dots,p_K\in [0,1]$.
Then $a_F(\epsilon)=\epsilon/K$, $b_F(\epsilon)=1-(1-\epsilon)^{1/K}$ and $c_F(\epsilon)=\epsilon$ for $\epsilon\in (0,1)$.

\end{proposition}

\subsection{The  averaging methods}
The aforementioned averaging methods of \cite{VW19} use the combining functions given by
\begin{align*}
    M_{r,K}(p_{1},\dots, p_{K})=\left(\frac{p_{1}^{r}+\dots+p_{K}^{r}}{K}\right)^{\frac{1}{r}},
\end{align*}
for $r\in\mathbb{R}\setminus\{0\}$, together  with its limit cases
\begin{align*}
   M_{-\infty,K}(p_{1},\dots, p_{K})&= \min\{p_{1},\dots,p_{K}\};\\
   M_{0,K}(p_{1},\dots, p_{K})&= \left(\prod_{i=1}^{K}p_{i}\right)^{\frac{1}{K}};\\
   M_{\infty,K}(p_{1},\dots, p_{K})&= \max\{p_{1},\dots,p_{K}\}.
\end{align*}
Some special cases of the combining functions  above are  $r=-\infty$ (minimum), $r=-1$ (harmonic mean), $r=0$ (geometric mean), $r=1$ (arithmetic mean) and $r=\infty$ (maximum); the cases $r\in\{-1,0,1\}$ are known as Platonic means. Note that $ M_{-\infty,K}$ gives rise to   the Bonferroni method, and the geometric mean  yields Fisher's method (\cite{fs1938methods}) under the independence assumption.
The harmonic mean p-value of \cite{W19} is a VSD method using the harmonic mean.

Since the mean function $M_{r,K}$ is homogeneous,
by Proposition \ref{p1},   the VAD threshold is  a linear function $a_{F}(x)=a_{r}x$, $x\in(0,1)$ for some $a_r>0$. The multipliers $a_{r}$ have been well studied in \cite{VW19}, and here we mainly focus on the cases of Platonic means and the Bonferroni method.
It is   known that {$a_{-\infty}=1/K$ and $a_{1}=1/2$.} For $r=0$ or $r=-1$, the values of $a_r$ and their asymptotic formulas  are calculated by  Propositions 4 and 6 of  \cite{VW19},
summarized below for $K\geq 3$.
\begin{enumerate}[(i)]
 \item
 For $F=M_{0,K}$,
 \begin{equation}\label{a0}
 a_{F}(x)=a_{0}x=c_{K}\exp\left(\frac{K-1}{1-Kc_{K}}\right)\times x,~~~x\in(0,1),
 \end{equation}
   where $c_{K}$ is the unique solution to the equation: $\log(1/c-(K-1))=K-K^{2}c$ for $c\in (0,1/K)$. Moreover, $a_{0}\geq 1/e$, and $a_{0}\rightarrow 1/e$ as $K\rightarrow \infty$.
 \item 
  For $F=M_{-1,K}$,
 \begin{equation}\label{a-1}
 a_{F}(x)=a_{-1}x=\frac{(y_{K}+1)K}{(y_{K}+K)^{2}}\times x,~~~x\in(0,1),
 \end{equation}
  where $y_{K}$ is the unique solution to the equation: $y^{2}=K((y+1)\log(y+1)-y)$ for $y\in (0,\infty)$. Moreover, $a_{-1}\geq (e\log {K})^{-1}$, and $a_{-1}\log K\rightarrow 1$ as $K\rightarrow \infty$.
\end{enumerate}

 To determine the VC threshold, it is easy to check that $c_{M_{r,K}}(x)=x$, $x\in(0,1)$ for all $r\in [-\infty,\infty]$, because the generalized mean of identical objects is equal to themselves; this obviously holds for all functions in the family of $M_{\phi,K}$.

Next,  we study  $b_r:= b_{M_{r,K}}$ or its approximate  form.
For this, we will use stable distributions   (e.g., {\cite{uchaikin2011chance} and} \cite{samoradnitsky2017stable}) below.
Let $F_\alpha$ be the stable distribution
with stability parameter $\alpha \in (0,2)$,
 skewness parameter $\beta=1$, scale parameter $\sigma=1$ and shift parameter $\mu=0$.
The characteristic function of $ F_{\alpha}$ is given by,  for $\theta \in \R$,
\begin{equation*}
    \int \exp(i\theta x)\d F_\alpha(x) =\begin{cases}
    \exp\left(-|\theta|^\alpha (1-i\sign(\theta)\tan\frac{\pi\alpha}{2})\right) \mbox{~~~~if $\alpha\neq 1$},\\
    \exp\left(-|\theta|(1+i\frac{2}{\pi}\sign(\theta)\log|\theta|)\right) \mbox{~~~~if $\alpha= 1$},
    \end{cases}
\end{equation*}
where $\sign(\cdot) $ is the sign function.
For $\alpha \ge 2$, let $F_\alpha$ stand for the standard normal distribution.

\begin{proposition}\label{p3}
Let  $b_{r}$ be the VI threshold of $M_{r,K}$, $r\in \R$.
\begin{enumerate}[(i)]
\item\label{p3i}
If $r<0$, then for $K\in \mathbb{N}_+$
\begin{equation}\label{VI:r<0}
 b_{r}  (\epsilon)\sim  K^{-1-1/r}\epsilon, \mbox{~~~~as $\epsilon \downarrow 0$}, 
\end{equation}
and for  $\epsilon\in (0,1)$,
$$  b_{r}  (\epsilon)\sim  \left(\left(C_{\alpha}F_{\alpha}^{-1}(1-\epsilon)+b_{K}\right)/K\right)^{\frac{1}{r}}, \mbox{~~~~as $K \rightarrow \infty$},$$
where $\alpha=-1/r>0$ and the constants $C_{\alpha}$ and $b_{K}$ are given in Table \ref{t1}.
\item
If $r=0$, then 
\begin{equation}\label{VI:r=0}
b_{r}(\epsilon)=\exp\left(-\frac{1}{2K}q_{1-\epsilon} \left(\chi^2_{2K}\right)\right).
\end{equation}
\item
If $r>0$, then for $K\in\mathbb{N}_+$,
$$  b_{r}  (\epsilon)= \frac{(\Gamma(1+K/p))^{1/K}\epsilon^{1/K}}{K^{1/r}\Gamma(1+1/p)}, \mbox{~~~~if $\epsilon \leq \frac{(\Gamma(1+1/p))^K}{\Gamma(1+K/p)}$},$$
where $\Gamma$ is the Gamma function.
For $\epsilon\in (0,1)$,
$$  b_{r}  (\epsilon)\sim  \left(\frac{\sigma}{\sqrt{K}}\Phi^{-1}(\epsilon)+\mu\right)^{\frac{1}{r}}, \mbox{~~~~as $K \rightarrow \infty$},$$
 where $ \mu=(r+1)^{-1}$ and $\sigma^{2}=r^2(1+2r)^{-1}(1+r)^{-2}.$
\end{enumerate}
\end{proposition}

 \begin{table}[htbp]
    \renewcommand\arraystretch{0.75}
    \caption{ Coefficients $C_{\alpha}$ and $b_{K}$ for $r=-1/\alpha <0$.}
\label{t1}
\centering
\footnotesize
\begin{tabular}{ccc}
 \hline
$r=-1/\alpha$&$C_{\alpha}$& $b_{K}$   \\
  \hline
  $-\frac{1}{2}< r<0$& $\left(K\left(\frac{\alpha}{\alpha-2}-\left(\frac{\alpha}{\alpha-1}\right)^{2}\right)\right)^{1/2}$ & $K\alpha/(\alpha-1)$   \\

  {$ r=-\frac{1}{2}$}& {$\sqrt{K\log K}$} & {$K\alpha/(\alpha-1)$}   \\

  $-1< r<-\frac{1}{2}$ & $K^{1/\alpha}\left(\Gamma(1-\alpha)\cos(\pi\alpha/2)\right)^{1/\alpha}$ & $K\alpha/(\alpha-1)$   \\

     $r=-1$ &  $K\pi/2$  & $\displaystyle \frac{\pi K^{2}}{2}\int_{1}^{\infty}\sin\left(\frac{2 x}{K\pi}\right) \alpha x^{-\alpha-1} \d x$  \\

     $r<-1$ & $K^{1/\alpha}\left(\Gamma(1-\alpha)\cos(\pi\alpha/2)\right)^{1/\alpha}$ & $0$   \\
   \hline
\end{tabular}
\end{table}

\subsection{The Cauchy combination method}

The Cauchy combination method is recently proposed by \cite{LX19} which relies on a special case of the generalized mean via $\phi=\mathcal C^{-1}$,
where $\mathcal C$ is the standard Cauchy cdf, that is,
\begin{align*}
    \mathcal C(x)&=\frac{1}{\pi}\arctan(x)+\frac{1}{2},~x\in \mathbb{R}; ~~~~
    \mathcal C^{-1}(p)=\tan\left(\pi\left(p-\frac{1}{2}\right)\right),~p\in(0,1).
\end{align*}
We denote this combining function by $M_{ \mathcal C,K}$  (instead of $M_{\mathcal C^{-1},K}$  for simplicity), namely,
\begin{align*}
    M_{ \mathcal C,K}(p_{1},\dots,p_{K}):=\mathcal C\left(\frac{1}{K}\sum_{i=1}^{K}\mathcal C^{-1}\left(p_{i}\right)\right).
\end{align*}
It is well known that the arithmetic average of either independent or comonotonic standard Cauchy random variables
follows again the standard Cauchy distribution.
This   feature allows the use of such a combination method to combine p-values
under uncertain dependence assumptions.
In addition, \cite{LX19} showed that under a bivariate normality assumption of the individual test statistics (i.e., a normal copula), the combined p-value has the same asymptotic behaviour as the one under the assumption of independence (see Theorem \ref{theorem2} (ii) below).

Since $\frac{1}{K}\sum_{i=1}^{K}\mathcal C^{-1}(U_{i})$ follows a standard Cauchy distribution if $U_{1},\dots,U_{K}\in \mathcal U$ are either independent or comonotonic,   we have $b_{F}(x)=c_{F}(x)=x$ for all $x\in(0,1)$.
This convenient feature will be studied in more details   in Section \ref{s4}.

By Definition \ref{def:1}, we  get, for $F=M_{\mathcal C,K}$,
\begin{align}\label{Cauchy}a_{F}(\epsilon)=\mathcal C \left(\inf\left\{q_{\epsilon}\left(\frac{1}{K}\sum_{i=1}^{K} \mathcal C^{-1}(U_{i})\right)\mid U_1,\dots,U_K\in \mathcal U\right\}\right).
\end{align}
The function $a_F$ does not admit an explicit formula, but it can be calculated via
 results from  robust risk aggregation (Corollary 3.7 in  \cite{WPY13}) as in the following proposition.
\begin{proposition}\label{p4} For $\epsilon \in (0,1/2)$, we have
\begin{equation}\label{VAD:cauchy}
a_F(\epsilon)=\mathcal C\left(- H_\epsilon (x_K)/K\right),
\end{equation}
where
$
H_\epsilon(x) = (K-1)\mathcal C^{-1}(1-\epsilon+(K-1)x)+\mathcal C^{-1}(1-x)$, $x\in (0,\epsilon/K)$, and
$x_K$ is the unique solution $x\in (0,\epsilon/K)$ to the equation
\begin{align*}
& K \int_{x}^{\epsilon/K}  H_\epsilon(t)\d t =(\epsilon- Kx)H(x).
\end{align*}
\end{proposition}

\subsection{The Simes method}
The     method  of \cite{S86} uses  the Simes function $ S_K$ in the order statistics family, given by $
  S_K(p_{1},\dots, p_{K}) =\min_{i\in\{1,\dots,K\}}\frac{K}{i}p_{(i)}.
 $
For $ F=S_K$, the results in \cite{hommel1983tests} together with Proposition \ref{p1} suggest that
$a_{F}(x)=x/\ell_K$ for $x\in(0,1)$. For independent  p-variables $P_1,\dots,P_K\in \mathcal U$, \cite{S86} obtained \begin{align*}
    \mathbb{P}\left(\min_{i\in\{1,\dots,K\}}\frac{K}{i}P_{(i)}>\epsilon\right)=1-\epsilon,~~~\epsilon \in(0,1),
    \end{align*}
which gives $b_{F}(x)=x$ for $x\in(0,1)$. For comonotonic p-variables $P_1,\dots,P_K\in \mathcal U$, it is clear that $S_K(P_1,\dots,P_K)= P_{(K)}$, which follows a standard uniform distribution, and hence we again have $c_{F}(x)=x$ for $x\in(0,1)$.
The validity of the Simes function using the VI (VC) threshold (called the Simes inequality) holds under many positive dependence structures; see e.g., \cite{sarkar1998some, sarkar2008simes}.

In the context of testing multiple hypotheses,
if  p-variables for several hypotheses are independent, the Benjamini-Hochberg procedure for controlling the false discovery rate (FDR) (\cite{benjamini1995controlling}) also relies on the Simes  function (in case all hypotheses are null).  Although the  Benjamini-Hochberg procedure is valid  for many practical models, to control the FDR under arbitrary dependence structure of p-variables, one needs to multiply the p-values by $\ell_{K} $, resulting in the Benjamini-Yekutieli procedure (\cite{benjamini2001control}). This  constant is exactly  $x/a_F(x)$,
and the function $a_F$ is called a reshaping function by \cite{RBWJ19} in the FDR context.

\section{Independence-comonotonicity balance}\label{s4}

As we have seen above, the Cauchy function and the Simes function both satisfy
 $b_F=c_F$, and hence
 the corresponding merging methods  are invariant under independence or comonotonicity assumption, an arguably convenient feature.
 Inspired by this observation,
  we introduce the property of \emph{independence-comonotonicity balance} for combining functions in this section.
This property distinguishes the Cauchy combination method and the Simes method from their corresponding classes $M_{\phi,K}$ and $S_{\alpha,K}$, respectively.

A combining function is said to be balanced between two different dependence structures of p-variables if the combined random
variable under the two dependence assumptions coincide in distribution. Recall that $U,V_{1},\dots,V_{K}$ are independent standard uniform random variables.
\begin{definition}
A combining function $F:[0,1]^{K}\rightarrow [0,\infty)$ is \emph{independence-comonotonicity balanced} (\emph{IC-balanced}) if $F(V_{1},\dots,V_{K})\laweq F(U,\dots,U).$
\end{definition}
As the VI and VC thresholds are the corresponding  quantile functions of $F(P_{1},\dots, P_{K})$,
we immediately conclude that a combining function $F:[0,1]^{K}\rightarrow [0,\infty)$ is IC-balanced if and only if $b_{F} =c_{F} $  on $(0,1]$; recall that $c_F$ is the identity for all functions in Section \ref{s3}.

IC-balanced methods have the same threshold  $b_{F}=c_{F}$ if the dependence structure of p-variables is a mixture of independence and comonotonicity, i.e., with the   copula
\begin{align}\label{copula}\lambda\prod_{i=1}^{n}x_i+(1-\lambda)\min_{i=1,\dots,n} x_i,~~~(x_1,\dots,x_n)\in [0,1]^n,
\end{align}
where  $\lambda\in [0,1]$. This is because $\p(F(U_1,\dots,U_K)\le b_F(\epsilon))$ is linear in the distribution of $(U_1,\dots,U_K)$.

For any combining function $F$, VI and VC thresholds generally yield more power to the test compared with the corresponding VAD threshold, but the gain of power may come with the invalidity due to model misspecification. If a combining function $F$ is IC-balanced, the validity is preserved under independence, comonotonicity  {and their mixtures}, and we may expect (without mathematical justification) that, to some extent, the size of the test can be controlled properly even if mild model misspecification exists. Therefore, the notion of IC-balance can be  interpreted as insensitivity to some specific type of model misspecification (e.g., dependence structure given in \eqref{copula})  for  VSD merging methods.

We have already seen in Section \ref{s3} that the Cauchy combination method and the Simes method are IC-balanced.
Below we show that they are the only IC-balanced methods among the two classes of combining functions based on generalized mean and order statistics.
\begin{theorem}\label{Proic-balance}
For a generalized mean function $M_{\phi,K}$ and
an order statistics function $S_{\alpha,K}$,
\begin{enumerate}[(i)]
\item    $M_{\phi,K}$ is IC-balanced for all $K\in\N$ if and only if it is the Cauchy combining function, i.e., $\phi(p) $ is a linear transform of $\tan\left(\pi\left(p-\frac{1}{2}\right)\right)$, $p\in(0,1)$;
\item $S_{\alpha,K}$ is IC-balanced if and only if it is a positive constant times the Simes function.
\end{enumerate}
\end{theorem}

The IC-balance of $M_{\phi,K}$ for some fixed $K$ (instead of all $K\in \N$) does not imply that $\phi $ is the quantile function of a Cauchy distribution; see the counter-example (Example \ref{ex:1}) in the supplementary material.
As a direct consequence of Theorem \ref{Proic-balance},  if $S_{\alpha,K}$ is IC-balanced, then $S_{\alpha,k}$ for $k=2,\dots, K-1$, are also IC-balanced (here we use the first $k$ components of $\alpha$); a similar statement does not hold in general for the generalized mean functions, also shown by Example \ref{ex:1}.
\begin{remark}
The property of IC-balance should be seen as a necessary but
not sufficient condition for a merging method to be insensitive to dependence between
independence and comonotonicity.  As shown by \cite{sarkar1998some},
the Simes method is valid for positive regression dependence, which is a large spectrum of dependence structures connecting independence and comonotonicity (larger than \eqref{copula}); on the other hand, the Cauchy {combination} method using VI threshold is valid under a bivariate Gaussian assumption   asymptotically  but not precisely (\cite{LX19}); see Theorem \ref{theorem2}  below and the simulation studies in Section \ref{sec:6}.
Instead of arguing for the practical usefulness of IC-balance,
we emphasize it as a necessary condition for insensitivity to dependence.
The main aim of Theorem \ref{Proic-balance} is, via this necessary condition, to pin down the unique role of the Simes and the Cauchy  {combination} methods among their respective generalized classes, thus justifying their advantages with respect to dependence.
\end{remark}

\section{Connecting the Simes, {the} harmonic {averaging} and  {the} Cauchy  {combination} methods}
\label{sec:4p}

As we have seen from Theorem \ref{Proic-balance},
the Cauchy and Simes combining functions are the only IC-balanced ones among the two classes considered in Section \ref{s3}.
Although the harmonic combining function does not satisfy $b_F=c_F$,  we observe empirically that  the harmonic averaging method and the Cauchy combination method
report very similar results in all simulations; see Section \ref{sec:6}.

 In this section, we explore the relationship among the three methods based on $S_K$, $M_{-1,K}$ and $M_{\mathcal C,K}$.
 We first show that the harmonic averaging method is equivalent to the Cauchy combination method asymptotically in a few senses.
Second,   we show the Simes function $S_K$ and the harmonic averaging function $M_{-1,K}$ are closely connected via
$M_{-1,K}\le S_K\le \ell_K M_{-1,K}$, where $\ell_K$ is given in \eqref{eq:ellK}.
Throughout this section, for fixed $K\in \N$, we write $a_{\mathcal C}= a_{M_{\mathcal C,K}}$, $a_{\mathcal S}= a_{S_{K}}$,
$a_{\mathcal H}= a_{M_{-1,K}}$ and similarly for $b_{\mathcal C}$, $b_{\mathcal S}$ and $b_{\mathcal H}$.

We will use  the following   assumption on the p-variables $U_{1},\dots,U_{K}\in \mathcal U$.
\begin{itemize}
    \item [(G)]
    For each $1\leq i<j\leq K$, $(U_{i},U_{j})$ follows a bivariate Gaussian copula (which can be different for each pair).
\end{itemize}
 The assumption  (G) is mild and is imposed by  \citet[Condition C.1]{LX19}. Note that condition (G) includes  independence and comonotonicity  as special cases. The following theorem confirms the close relationship between the harmonic averaging method  and the Cauchy combination method.
 Recall that   the VC thresholds for both methods are the identity function,
 and thus it suffices to look at VAD  and VI thresholds.
\begin{theorem}\label{theorem2}
For fixed $K\in \N$, the harmonic  {averaging} and the Cauchy combination methods are asymptotically equivalent in the following senses:
\begin{enumerate}[(i)]
    \item  If $\min_{i\in\{1,\dots,K\}}p_{i}\downarrow 0$ and $\max_{i\in\{1,\dots,K\}}p_{i}\leq c$ for some fixed $c\in (0,1)$, then
    $$\frac{M_{\mathcal C,K}(p_{1},\dots,p_{K})}{M_{-1,K}(p_{1},\dots,p_{K})}\rightarrow 1.$$
      \item
For $K$ standard uniform random variables $U_{1},\dots,U_{K}$ satisfying condition (G),
\begin{equation}
\label{e11}
 \mathbb{P}\left(M_{\mathcal C,K}(U_{1},\dots,U_{K})<\epsilon\right) \sim \mathbb{P}\left(M_{-1,K}(U_{1},\dots,U_{K})<\epsilon\right) \sim \epsilon,~\text{as}~\epsilon\downarrow 0.
\end{equation}
    In particular, $b_{\mathcal C} (\epsilon) \sim b_{\mathcal H} (\epsilon)$ as $\epsilon \downarrow 0$.
    \item   $a_{\mathcal C}(\epsilon)\sim a_{\mathcal H}(\epsilon)$ as   $\epsilon\downarrow 0$.
\item For $r\neq -1$,
     $$\frac{M_{\mathcal C,K}(p_{1},\dots,p_{K})}{M_{r,K}(p_{1},\dots,p_{K})}\not\rightarrow 1,~ \text{as}\max_{i\in\{1,\dots,K\}}p_{i}\downarrow 0.$$
\end{enumerate}
\end{theorem}
\begin{remark}
The statement $ \mathbb{P}\left(M_{\mathcal C,K}(U_{1},\dots,U_{K})<\epsilon\right) \sim \epsilon $
in Theorem \ref{theorem2} (ii)  is implied by Theorem 1 of \cite{LX19}, which gives  the same convergence  rate for the weighted Cauchy combination method.
For the weighted harmonic averaging method, we have a similar result  (see \eqref{m-1} in the supplementary material):
For standard uniform random variables $U_{1},\dots,U_{K}$ satisfying condition (G) and any $(w_1,\dots,w_K)\in [0,1]^K$ with $\sum_{i=1}^K w_i =1$, we have
$${
 \mathbb{P}\left( \sum_{i=1}^K w_i U_i^{-1} >1/\epsilon\right) \sim \epsilon,~\text{as}~\epsilon\downarrow 0.}
$$
We omit a discussion on weighted merging methods as the focus of this paper is comparing symmetric combination functions.
\end{remark}

The first statement of Theorem \ref{theorem2} means that, if at least one of realized p-values are close to 0, the harmonic averaging and the Cauchy combining functions will produce very close numerical results. This case is likely to happen in high-dimensional situations where the number of p-variables is very large. As the condition (G) for (ii) in   Theorem \ref{theorem2} is arguably mild,  the thresholds of the two methods are similar for a small significance level under a wide range of dependence structures of p-variables (including independence and comonotonicity).
Therefore, if the significance level is small, one likely arrives at the same statistical conclusions on the hypothesis testing by using  either method.
The third result in Theorem \ref{theorem2} illustrates the equivalence between the VAD thresholds of the harmonic averaging and the Cauchy combination methods as the significance level goes to 0.
The final result in Theorem \ref{theorem2} shows that among all averaging methods, the harmonic averaging method is the only one that is asymptotically equivalent to  the Cauchy combination method.
\begin{remark}
 We note that the   equivalence $$ \mathbb{P}\left(M_{\mathcal C,K}(U_{1},\dots,U_{K})<\epsilon\right) \sim \mathbb{P}\left(M_{-1,K}(U_{1},\dots,U_{K})<\epsilon\right)$$ in \eqref{e11}   does not always hold under arbitrary dependence structures.
Since the Cauchy distribution is symmetric, it is possible that
$\p (\mathcal C^{-1}(U_1)+\dots + \mathcal C^{-1}(U_K)=0)=1$
for some $U_1,\dots,U_K\in \mathcal U$, implying $\p(M_{\mathcal C,K}(U_1,\dots,U_K)<1/2)=0$.
Indeed,  Theorem 4.2 of \cite{puccetti2019centers} implies that there exist $K$ standard Cauchy random variables whose sum is a constant $c$,  for each  $c\in [-K\log(K-1)/\pi, K\log(K-1)/\pi]$.
On the other hand,  $\p( M_{-1,K}(U_1,\dots,U_K)<\epsilon)>0$
for all $\epsilon>0$ and all $U_1,\dots,U_K\in \mathcal U$.  Thus, $ \mathbb{P}\left(M_{\mathcal C,K}(U_{1},\dots,U_{K})<\epsilon\right) \sim \mathbb{P}\left(M_{-1,K}(U_{1},\dots,U_{K})<\epsilon\right)$ does not hold.
\end{remark}
\begin{remark}
The equivalence in Theorem \ref{theorem2} (ii) relies on the p-variables being uniform on $[0,1]$.
For p-variables that are stochastically larger than uniform, the behaviour of the Cauchy combination  method and that of the harmonic averaging method may diverge; nevertheless, by Theorem \ref{theorem2} (i), for a realized vector of p-values with at least one very small component, the two methods would produce similar values.
\end{remark}

The next result reveals an intimate relationship between the Simes and the harmonic averaging methods.
\begin{theorem}\label{sim_har}
For $p_{1},\dots,p_{K}\in[0,1]$, $$M_{-1,K}(p_{1},\dots,p_{K}) \le S_K(p_{1},\dots,p_{K}) \le \ell_K M_{-1,K}(p_{1},\dots,p_{K}).$$
 The first inequality holds as an equality if $p_{1}=\dots=p_{K}$.
The second inequality holds as an equality if $p_{1}=p_{k}/k$ for $k=2,\dots,K$.
As a consequence, $a_{\mathcal S}/a_{\mathcal H} \in [1,\ell_K]$ and $b_{\mathcal S} / b_{\mathcal H} \in [1,\ell_K]$.
\end{theorem}

 By Proposition \ref{p3} (i), the VI threshold of the harmonic averaging method
 satisfies $b_{\mathcal H}(\epsilon) \sim \epsilon= b_{\mathcal S}$ as $\epsilon\downarrow 0$.
Using Theorem \ref{sim_har}, we further know that $b_{\mathcal H}(\epsilon)<\epsilon$ (the inequality is strict since $M_{-1,K}<S_K$ has probability $1$ for independent p-variables).
Therefore, we cannot directly use the asymptotic VI threshold $\epsilon$ of the harmonic averaging method, which needs to be corrected; see \cite{W19}.

To summarize the results in this section, the Cauchy combining function and the harmonic averaging function are very similar in several senses, and the Simes function is more conservative  than the harmonic averaging function. Empirically, we see that the Simes function is only slightly more conservative; see Section \ref{sec:6}.

\section{Prices for validity}\label{s5}

For a given set of realized p-values, the decision to the hypothesis testing for some specific combining function will  be determined by the corresponding threshold. The VAD method can always control the size below the significance level; VSD methods may not have the correct size, but they yield more power than the VAD method. Therefore, there is always a trade-off between   validity and   efficiency, thus a price for validity.

For a combining function $F$ and {$K$ standard uniform random variables $U_{1},\dots,U_{K}$ with some specific dependence assumption (e.g., independence, comonotonicity, or condition (G)), let $g_{F}$ be the VSD threshold, i.e., $g_{F}(\epsilon)=q_{\epsilon}(F(U_{1},\dots,U_{K}) )$. Let $a_{F}$ be defined as in \eqref{e1p}.} For some fixed $\epsilon\in(0,1)$, the ratio  $ {g_{F}(\epsilon)}/{a_{F}(\epsilon)}$  is called the \emph{price for validity} under the corresponding dependence assumption of the p-variables.
For instance, $b_{F}(\epsilon)/a_{F}(\epsilon)$ is the price paid for validity under independence assumption and $c_{F}(\epsilon)/a_{F}(\epsilon)$ is the corresponding price under the comonotonicity assumption.
For a   specific application, one may consider the price for validity under other dependence assumptions. The calculation of the price for validity serves for two purposes:
\begin{enumerate}[i]
     \vspace{-0.3cm}\item (Power gain/loss): On the one hand, if additional information on the dependence structure of the p-values is available, the price for validity can be used as a measure for the gain of power   from the dependence information. On the other hand, if the dependence information is not available or credible, the price can be used to measure the power loss by  switching to the VAD threshold.
     \vspace{-0.3cm}\item (Sensitivity to model misspecification): If the dependence structure is ambiguous, VAD thresholds should be used. A small price for validity indicates
    that  a relatively small change of threshold due to the model ambiguity. Hence, the price for validity can be used as a tool to assess the sensitivity of VSD methods to model misspecification.
\end{enumerate}
\begin{remark}\label{remark6}
{Instead of using the price for validity,  a more direct way to assess the trade-off between using VSD and VAD methods is comparing the sizes}, i.e., $\mathbb{P}(F(P_{1},\dots,P_{K})<g_F(\epsilon))/\mathbb{P}(F(P_{1},\dots,P_{K})<a_F(\epsilon))$, where the dependence of p-variables $P_{1},\dots,P_{K}$ corresponds to the VSD method.  More precisely, for a fixed $\epsilon\in (0,1)$, the ratio of sizes is ${\epsilon}/{g_F^{-1}(a_F(\epsilon))}$, where $g_F^{-1}$ is the (generalized) inverse of $g_F$. The connection between the price for validity and the ratio of sizes is explained below.

\begin{enumerate}[(i)]
 \item For the Simes and the Cauchy combination methods, {the ratios of sizes under independence and  comonotonicity  are identical to the corresponding price for validity since $b_F$ and $c_F$ are identity functions.}
 \item For the averaging methods, {the ratios of sizes under comonotonicity are identical to the price for validity since $c_F$ is identity. The ratios of sizes under independence} may be different from $ {b_F(\epsilon)}/{a_F(\epsilon)}$; however, by letting $\delta=a_F(\epsilon)$, we have ($a_F$ is strictly increasing in all cases we consider)
$$
\frac{\epsilon}{b_F^{-1}(a_F(\epsilon))} = \frac{a_F^{-1}(\delta)}{b_F^{-1}(\delta)}.
$$
This is  very similar to $ {b_F(\epsilon)}/{a_F(\epsilon)}$; it is a matter of looking at the ratio of threshold functions or that of their inverses. {In fact, if $r<0$, by Proposition \ref{p3}, we have,
$$\frac{\epsilon}{b_F^{-1}(a_F(\epsilon))}\sim \frac{b_F(\epsilon)}{a_F(\epsilon)},~\epsilon\downarrow 0,$$
which suggests that the ratio of sizes is almost the same as the price for validity under independence for small significance levels.}
\end{enumerate}
\end{remark}

We use the Bonferroni method based on the combining function $F=M_{-\infty,K}$ as an example to illustrate the above idea.
Using Proposition \ref{prop:bonf} and noting that $K(1-(1-\epsilon)^{1/K})\sim \epsilon$ as $\epsilon \downarrow 0$,  we obtain that
 the prices for validity  of the Bonferroni method
satisfy
$c_F(\epsilon)/a_F(\epsilon) = K$ for $\epsilon \in (0,1 )$ and
$b_{F}(\epsilon)/a_{F}(\epsilon)\rightarrow 1$ as $\epsilon\downarrow 0$.
Therefore,
for a small $\epsilon$ close to 0, the price for validity under the independence assumption  is close to 1 while the price for validity under the comotonicity assumption   increases linearly as the number of p-variables increases.
This means a model misspecification of independence is not affecting the Bonferroni method much,
whereas a model misspecification of comonotonicity greatly affects the statistical conclusion of the Bonferroni method.

 Next we numerically calculate the prices for validity under independence and comonotonicity assumptions for various merging methods using results in Section \ref{s3}.
 We consider  the Bonferroni, the harmonic averaging, the geometric averaging, the Cauchy combination,  the Simes, and the negative-quartic
 (using $M_{-4,K}$, a compromise between Bonferroni and harmonic averaging)  methods. The (asymptotic) VAD and VI thresholds of these methods are summarized in Table \ref{threshold}. The VC threshold  is identity for all these methods. The VAD threshold of the negative-quartic method is given by Proposition 5 of \cite{VW20}.
  Numerical results on the prices for validity are reported  in Table \ref{t2} for $\epsilon =0.01$. Although some of the VAD thresholds in Table \ref{threshold} do not have explicit forms, the numerical computation is very fast.
  The results  for  $\epsilon=0.05$ and {$\epsilon=0.0001$} are similar and  reported in  {Tables \ref{t2p} and \ref{t2p1}} in the supplementary material.
\begin{table}[htbp]
  \caption{Thresholds for $K$ p-variables at significance level $\epsilon\in(0,1)$.}
\label{threshold}
\centering \small
\begin{tabular}{lcccccc}

   \hline
                    & {\scriptsize Bonferroni} & {\scriptsize Negative-quartic} & {\scriptsize Simes} & {\scriptsize Cauchy} & {\scriptsize Harmonic} & {\scriptsize Geometric}\\
   \hline
   $a_{F}(\epsilon)$  & $\epsilon/K$ & $\frac{3}{4}K^{-\frac{3}{4}}\epsilon$ & $\epsilon/\ell_K$ & \eqref{VAD:cauchy} & \eqref{a-1} & \eqref{a0} \\

   $b_{F}(\epsilon)$  & $1-(1-\epsilon)^{1/K}$ & \eqref{VI:r<0} & $\epsilon$ & $\epsilon$ & \eqref{VI:r<0} & \eqref{VI:r=0} \\
   \hline
   
\end{tabular}

\end{table}

 \begin{table}[htbp]
  \caption{ $b_{F}(\epsilon)/a_{F}(\epsilon)$ and $c_{F}(\epsilon)/a_{F}(\epsilon)$   for $\epsilon=0.01$ and $K\in \{50,100,200,400\}$}
\label{t2}
\centering \small
\begin{tabular}{lcccccccc}
\hline
& \multicolumn{2}{c}{$K=50$} & \multicolumn{2}{c}{$K=100$}& \multicolumn{2}{c}{$K=200$}& \multicolumn{2}{c}{$K=400$}\\
   \hline
                    & $b_{F}/a_{F}$ & $c_{F}/a_{F}$ & $b_{F}/a_{F}$ & $c_{F}/a_{F}$ & $b_{F}/a_{F}$ & $c_{F}/a_{F}$ & $b_{F}/a_{F}$ & $c_{F}/a_{F}$ \\

   {Bonferroni} & 1.005 & 50.000 & 1.005 & 100.000& 1.005 & 200.000& 1.005 & 400.000\\

    {Negative-quartic} & 1.340 & 25.071 & 1.340 & 42.164& 1.340 & 70.911 & 1.340 & 119.257 \\

      {Simes} & 4.499 & 4.499 & 5.187 & 5.187& 5.878 & 5.878 & 6.570 & 6.570 \\

    {Cauchy} & 6.625 & 6.625 & 7.465 & 7.465& 8.277 & 8.277  & 9.058 & 9.058\\

    {Harmonic} & 6.658 & 6.625 & 7.496 & 7.459& 8.314 & 8.273 & 9.117 & 9.072\\

    {Geometric} & 69.903 & 2.718 & 78.096 & 2.718& 84.214 & 2.718  & 88.694 & 2.718\\
   \hline
\end{tabular}
\end{table}

The Bonferroni and the negative-quartic methods
pay much lower price under the independence assumption than the comonotonicity assumption,
and the geometric averaging method is the absolute opposite.
On the other hand, the harmonic averaging, the Simes and the Cauchy combination methods have relatively small prices under both independence and comonotonicity assumptions and their prices increase at moderate rates as $K$ increases, compared to other methods. In particular, the harmonic averaging and the Cauchy combination methods have very similar performance (cf.~Theorem \ref{theorem2}) and their prices are slightly larger than that of the Simes method. If mild model misspecification exists, it may be safer to choose one of the harmonic averaging, the Simes and the Cauchy combination methods and use the corresponding VAD threshold without losing much power. {The prices for validity in Table \ref{t2} can also be interpreted as inflations of sizes by using VSD threshold against VAD threshold except the geometric averaging method (see Remark \ref{remark6}).}

 Next, we show that the prices for validity of the harmonic averaging, the Cauchy combination and the Simes methods behave like $\log K$ for $K$ large enough and $\epsilon$ small enough.
\begin{proposition}\label{prop:logK}
For $\epsilon \in (0,1)$, the prices for validity satisfy:
\begin{enumerate}[(i)]
    \item
    For the harmonic averaging method, $F=M_{-1,K}$,
    $$ \lim_{\delta\downarrow 0}\frac{b_{F}(\delta)}{a_{F}(\delta)}  =   \frac{c_{F}(\epsilon)}{a_{F}(\epsilon) } \sim \log K,~\text{as}~K\to\infty.$$
    \item
    For the Cauchy combination method, $F=M_{\mathcal C,K}$,
        $$ \lim_{\delta\downarrow 0}\frac{b_{F}(\delta)}{a_{F}(\delta)} =   \lim_{\delta\downarrow 0} \frac{c_{F}(\delta)}{a_{F}(\delta) } \sim \log K,~\text{as}~K\to\infty.$$
    \item
     For the Simes method, $F=S_{K}$,
             $$  \frac{b_{F}(\epsilon)}{a_{F}(\epsilon)}  = \frac{c_{F}(\epsilon)}{a_{F}(\epsilon) } \sim \log K,~\text{as}~K\to\infty.$$
\end{enumerate}
 \end{proposition}

 Numerical values of the ratios between the price  for validity under independence assumption  and $\log K$  are reported in Table \ref{t3}; the results for the corresponding ratios under comonotonicity assumption  are similar for these methods. The Simes method has the fastest convergence rate among the three methods. The ratios for the harmonic averaging and the Cauchy combination methods converge quite slowly and have similar rates.
 This fact can also be explained by Theorem \ref{sim_har}, where we see that the Simes function is generally larger than the harmonic averaging function.

Based on Proposition \ref{prop:logK}, one may be tempted to use $b_{F}/\log K$ as the corrected critical value under model misspecification; however, for the harmonic averaging and the Cauchy combination methods, the asymptotic rate of $\log K$ can only be expected for very large $K$ (instead, $1.7\log K$ works for $K\ge 100$).

  \begin{table}[htbp]
    \caption{Numerical values of $\frac{1}{\log(K)}\frac{b_{F}(\epsilon)}{a_{F}(\epsilon)}$ for the Simes, the Cauchy combination and the harmonic averaging methods.}
\label{t3}
\centering \small

\begin{tabular}{lccccccc}
 \hline
 &$\epsilon$& $K=10$ & 20 & 50 & 100 & 200 & 500 \\
  \hline
   \multirow{2}{*}{Simes}& $0.05$ & 1.272035 & 1.200955 & 1.150097 & 1.126425 & 1.109415 & 1.093041 \\
   &$0.01$  & 1.272035 & 1.200955 & 1.150097 & 1.126425& 1.109415 & 1.093041 \\

   \multirow{2}{*}{Cauchy} & $0.05$ & 1.979572 & 1.82826 & 1.693025 & 1.620527 & 1.561670 & 1.511264 \\
  & $0.01$  & 1.980144 & 1.828822 & 1.693562 & 1.621011 & 1.562121 & 1.504288 \\

    \multirow{2}{*}{Harmonic} & $0.05$  & 2.026308 & 1.873762 & 1.73641 & 1.661098 & 1.601539 & 1.539448 \\
 & $0.01$ & 1.989255 & 1.837605 & 1.701851 & 1.627702 & 1.569179 & 1.508248 \\
 \hline
\end{tabular}
\end{table}

\section{Simulations and a real data example}\label{sec:6}

\subsection{Simulation studies}
We conduct $K$ one-sided z-tests of the null hypothesis: $\mu_{i}=0$ against the alternative hypothesis $\mu_{i}>0$, $i=1,\dots,K$, using the test statistic $X_{i}$ and the p-value $p_{i}$ from the \emph{i}th test, $i=1,\dots,K$. The tests are formulated as the following: $$p_{i}=\Phi(X_{i}),~~X_{i}=\rho Z+\sqrt{1-\rho^2}Z_{i}-\mu_{i},~~i=1,\dots,K.$$
where $\Phi$ is the  standard normal distribution function, $Z,Z_{1},\dots,Z_{K}$ are iid standard normal random variables,   $\mu_{i}\geq 0$, $i=1,\dots,K$, and $\rho$ is a parameter in $[0,1]$. Note that for $\rho=0$, the p-variables are independent, and $\rho=1$ corresponds to the case where  p-variables are comonotonic.

Let $K\in\{50,200\}$ and set the significance level $\epsilon=0.01$. To see how different dependence structures and signals   affect the size and the power for various methods using both VAD and VSD thresholds, the rejection probabilities (RPs) are computed over $\rho\in[0,1]$ under the following four cases:
\begin{enumerate}[(i)]
 \vspace{-0.3cm}\item (no signal)
  $100\%$ of $\mu_{i}$'s  are 0;
 \vspace{-0.3cm}\item (needle in a haystack)
 $98\%$ of $\mu_{i}$'s are 0 and $2\%$ of $\mu_{i}$'s are 4;

 \vspace{-0.3cm}\item (sparse signal)
$90\%$ of  $\mu_{i}$'s are 0 and 10\% of $\mu_{i}$'s are 3;
 \vspace{-0.3cm}\item (dense signal)
$100\%$ of $\mu_{i}$'s  are 2.
\end{enumerate}

\begin{figure}[htbp]
\centering
\includegraphics[height=3.71cm, trim={0 40 0 20},clip]{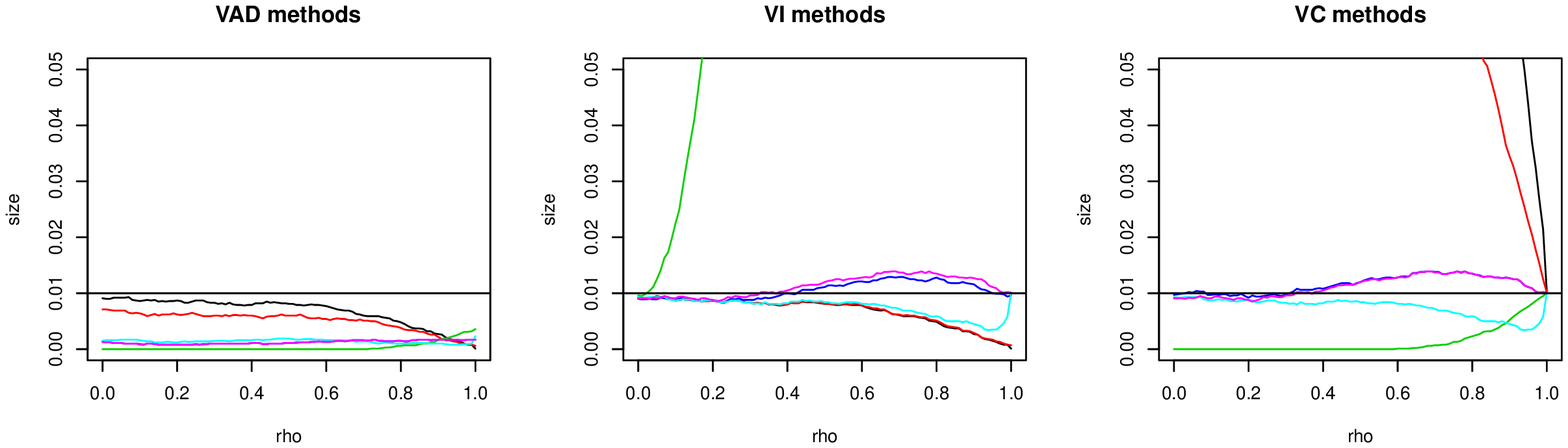}\\
\includegraphics[height=4.1cm, trim={0 0 0 40},clip]{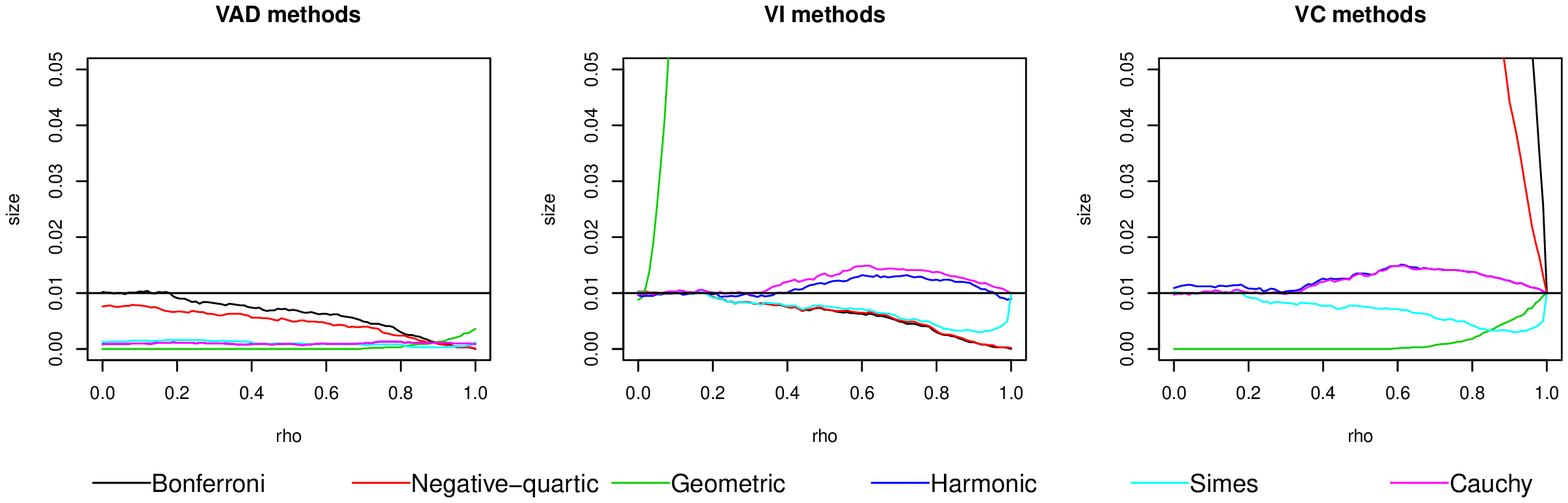}
\caption{Case (i): size (top: $K=50$, bottom: $K=200$)}
\label{f1}
\end{figure}

\begin{figure}[htbp]
\centering
\includegraphics[height=3.71cm, trim={0 40 0 20},clip]{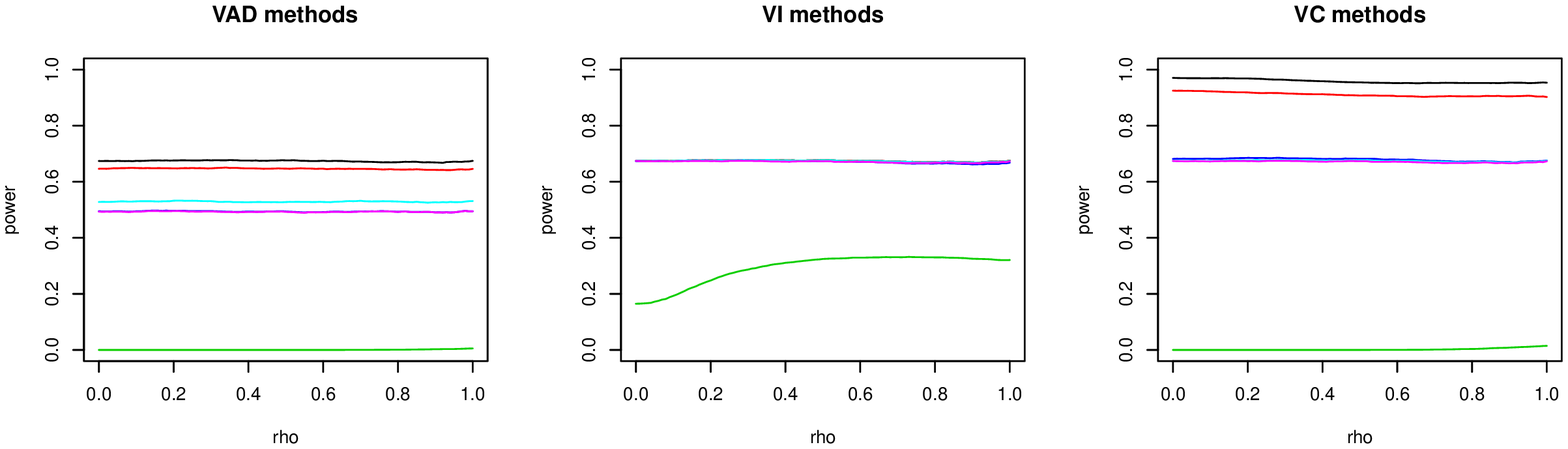}\\
\includegraphics[height=4.1cm, trim={0 0 0 40},clip]{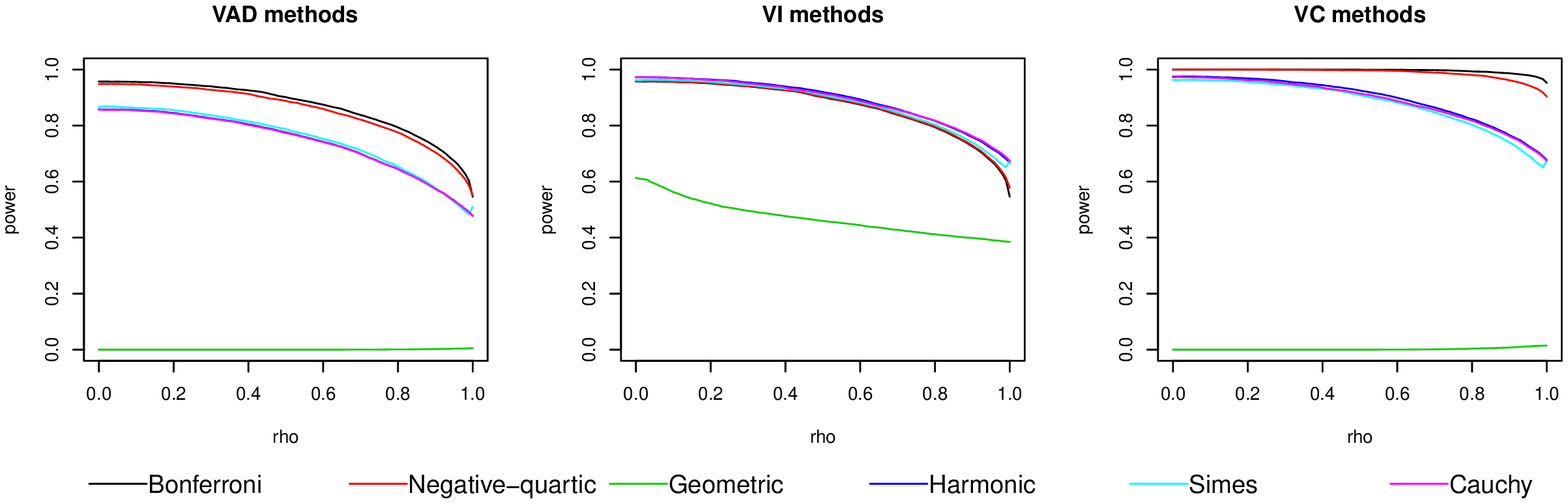}
\caption{Case (ii): needle in a haystack (top: $K=50$, bottom: $K=200$)}
\label{f2}
\end{figure}

\begin{figure}[htbp]
\centering
\includegraphics[height=3.71cm, trim={0 40 0 20},clip]{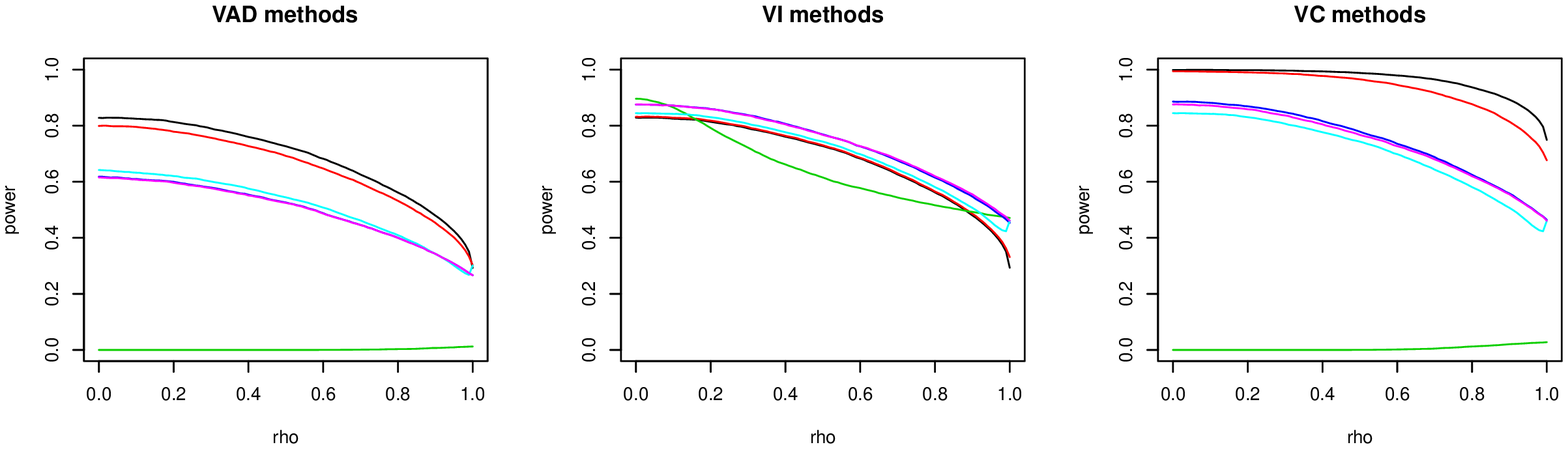}\\
\includegraphics[height=4.1cm, trim={0 0 0 40},clip]{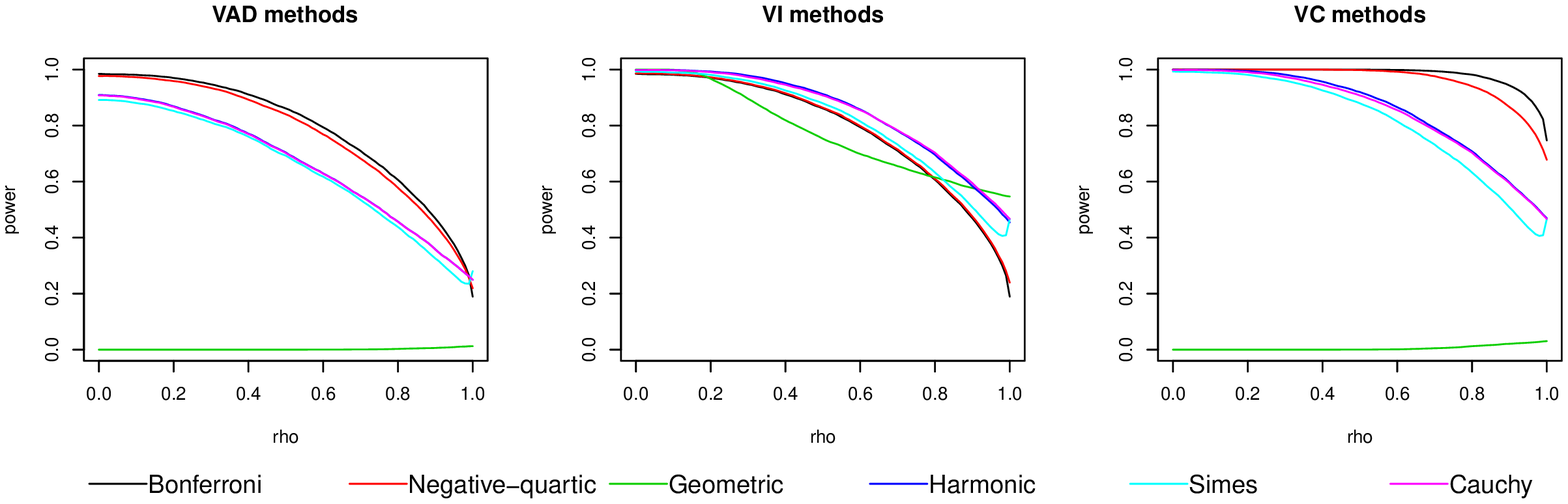}
\caption{Case (iii): sparse signal (top: $K=50$, bottom: $K=200$)}
\label{f3}

\end{figure}

\begin{figure}[htbp]
\centering
\includegraphics[height=3.71cm, trim={0 40 0 20},clip]{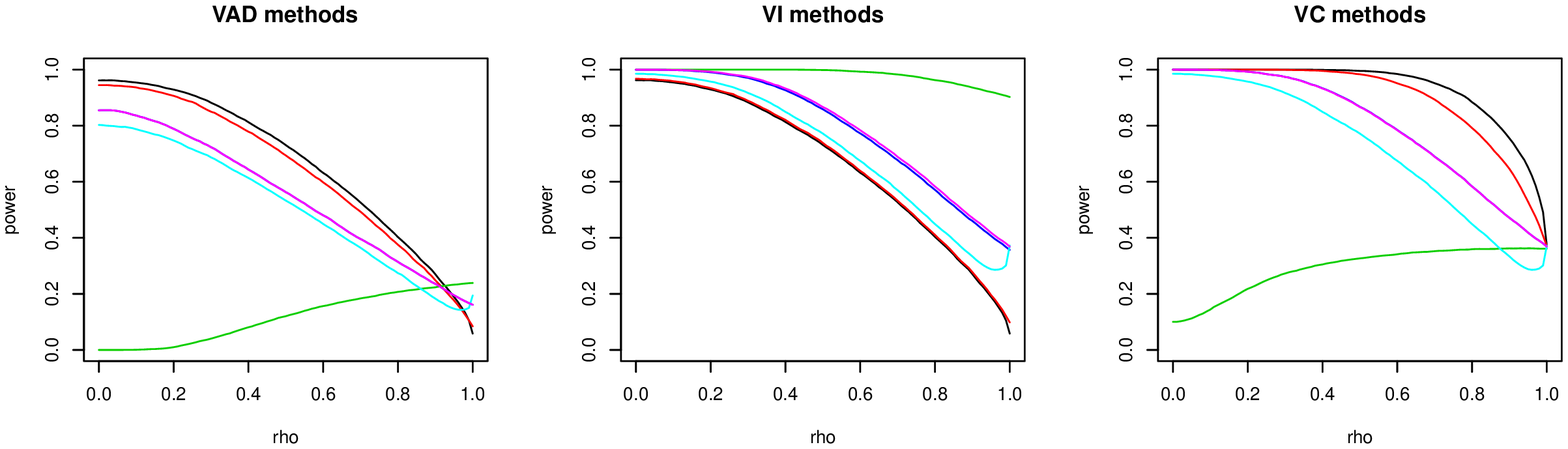}\\
\includegraphics[height=4.1cm, trim={0 0 0 40},clip]{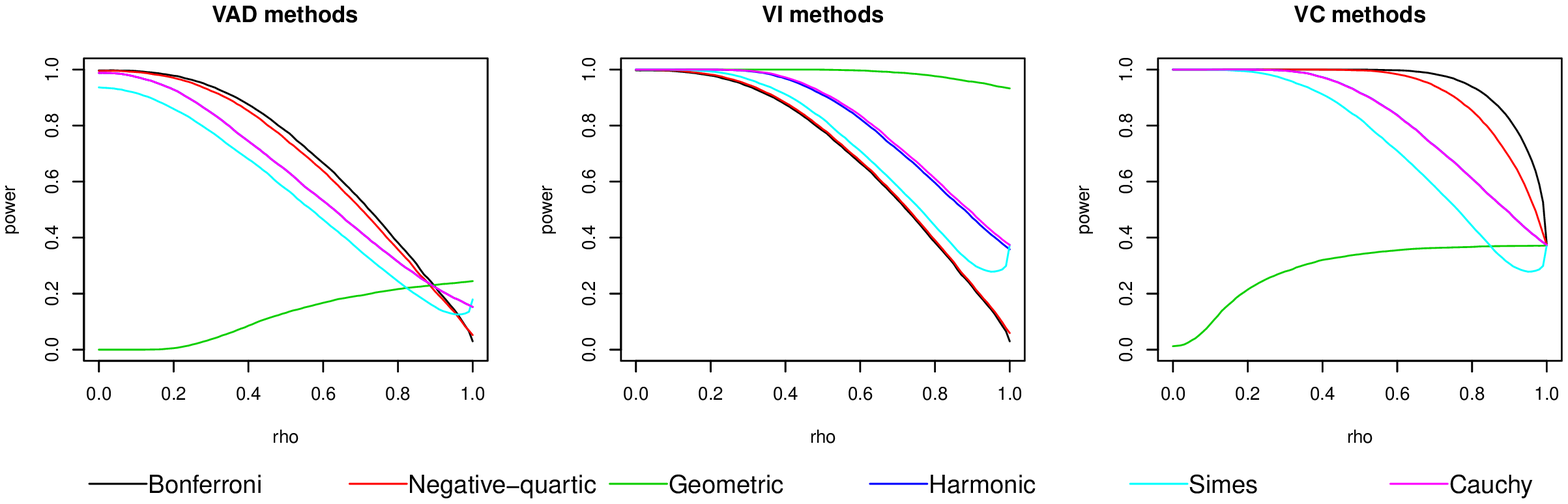}
\caption{Case (iv): dense signal (top: $K=50$, bottom: $K=200$)}
\label{f4}
\end{figure}

 The RP corresponds to the size under case (i), and it corresponds to the power under (ii), (iii) and (iv). The RP is computed as the ratio between the number of the combined values which are less than the critical threshold and the number of simulations for some $\rho\in[0,1]$, that is,
$$\mathrm{RP}=\frac{\sum_{i=1}^{N}\id_{\{F_{i}<g(\epsilon)\}}}{N},$$
where $N$ is the number of simulations and is equal to 15000 in our study,   $F_{i}$ is the realized value of the combining function for the $i$-th simulation, $i=1,\dots,N$, and $g(\epsilon)$ is the corresponding critical value.
For $\rho\in[0,1]$, graphs of RPs for different combining methods are drawn using VAD thresholds and VSD thresholds.
Some observations from Figures \ref{f1}-\ref{f4}   are made below, and
 those on the averaging methods using $M_{r,K}$ are consistent with the observations in \cite{VW19}.
\begin{enumerate}
 \item
All VAD methods give sizes less than $\epsilon=0.01$ as expected.
Using VAD thresholds,  the Bonferroni, the harmonic averaging, the Cauchy combination and the Simes methods have good powers.
\item
The Simes method using thresholds $b_F$ or $c_F$ reports the right size for all values of $\rho$. \cite{sarkar1998some}   showed the validity of the Simes method in the so-called MTP$_2$ class  including multivariate normal distributions with nonnegative correlations (the setting of our simulation).
\item
Using thresholds $b_F$ or $c_F$,  the harmonic averaging and Cauchy combination methods perform similarly with sizes possibly larger than $0.01$ (see Theorems \ref{theorem2} and \ref{sim_har}).
\item 
The geometric averaging method  using $b_F$
and the Bonferroni and negative-quartic methods using  $c_F$ do not yield correct sizes under model misspecification, and the sizes increase rapidly as the misspecification gets bigger.

\item 
Using $b_F$ or $c_F$, the harmonic averaging, the Cauchy combination and the Simes methods have good performances on capturing the signals.
\end{enumerate}

\subsection{Real data analysis}\label{rda}
We apply several merging methods to a genomewide study to compare their performances. We use the dataset of p-values of \cite{storey2003statistical} which contains 3170 p-values computed based on the data from \cite{hedenfalk2001gene} for testing whether genes are differentially expressed between BRCA1- and BRCA2-mutation-positive tumors. As mentioned in Section \ref{S2}, $g^{-1}\circ F(P_{1},\dots,P_{K})$ is a p-variable if the threshold $g$ is strictly increasing, and it is the quantity we choose to compare combined p-values for different methods.

\begin{figure}[htbp]
\centering
\includegraphics[width=15.1cm, trim={20 2 2 2},clip]{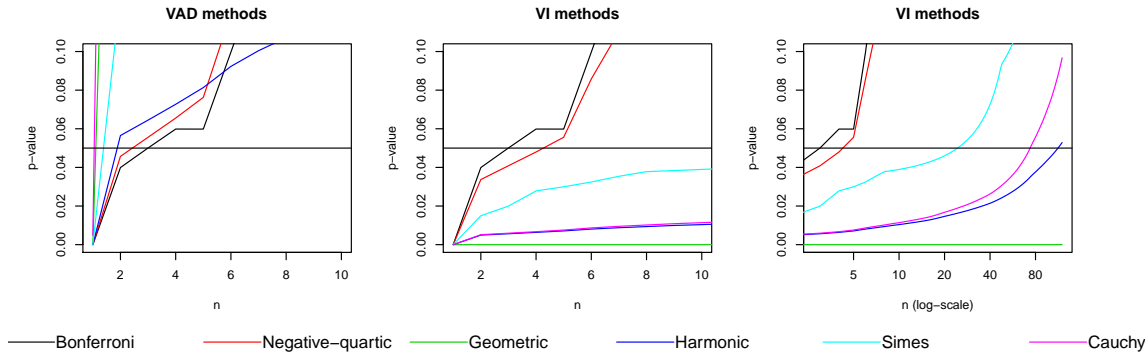}
\caption{Combined p-value after removing $n$  smallest p-values}
\label{f5}
\end{figure}

For each method, we calculate the combined p-value, and remove the smallest p-value from the dataset.  Repeat this procedure until the resulting combined p-value loses significance.
Using the Bonferroni combining function, this leads to the Bonferroni-Holm (BH) procedure (\cite{Holm:1979}); thus we  mimic the BH procedure for other methods in a naive manner.
The rough interpretation is to report the number of significant discoveries (this procedure generally does not control the family-wise error rate (FWER); to control FWER  one needs to use a generalized BH procedure as in \cite{VW19} or \cite{goeman2019simultaneous}. This procedure can be seen as a lower confidence bound from a closed testing perspective).
For a visual comparison of detection power, the combined p-values against the numbers of removed p-values are plotted in Figure \ref{f5}, where we use both the VAD and the VI thresholds (comonotonicity is obviously unrealistic here). {In the third panel of Figure 5, we present  the number of omitted p-values in log-scale for better visualization.}

All VAD methods lose significance at $\epsilon=0.05$ after omitting the first or the second smallest p-value (the smallest p-value is $0$ and the second smallest is $1.26\times 10^{-5}$). Using thresholds $b_F$ for independence, the Bonferroni and the negative quartic methods behave similarly to their VAD versions (as their price for validity is close to $1$). In contrast, the Simes, the Cauchy combination and the harmonic averaging methods lose significance at $\epsilon=0.05$ after removing around 20, 70 and 110 p-values respectively. The geometric averaging method (Fisher's)  exceeds 0.05 only after removing around 400 p-values.
However, this method relies  heavily on the independence assumption, which is impossible to verify from just one set of p-values.

\section{Concluding remarks}\label{S6}
We discussed two  aspects of merging p-values:
the impact of the dependence structure on the critical thresholds and  the trade-off between validity and efficiency.
 The Cauchy combination method and the Simes method are shown to be the only IC-balanced members among the generalized mean class and the order statistics class of combining functions.
The harmonic averaging   and the Cauchy combination methods are asymptotically equivalent, and the Simes and the harmonic averaging methods have simple algebraic relationship.
For the above three methods, the prices for validity under independence (comonotonicity) assumption all behaves like $\log K$ for large $K$. Moreover, these methods lose moderate amount of power if VAD thresholds are used, and their performance against model misspecification is better than other methods.  This explains the wide applications of these methods in different statistical procedures.

Merging p-values is not only useful for testing a single hypothesis, but also important in testing multiple hypotheses, controlling false discovery rate (\cite{benjamini1995controlling}, \cite{benjamini2001control}),  and exploratory research (\cite{goeman2011multiple}, \cite{goeman2019simultaneous}). 
In many situations especially involving a large number of hypotheses and tests, dependence information is hardly available.
The results in our paper offer some insights, especially in terms of gain/loss of validity and power,  on how the absence of such information
influences different statistical procedures of merging p-values.

In many practical applications, p-values arrive sequentially in time, and the
  existence of the $n$-th p-variable may depend on previously observed p-values (only promising experiments may be continued); thus the number of experiments to combine is a stopping time. Unfortunately,
the current merging method of p-values discussed in this paper cannot be used to sequentially update p-values with  arbitrary stopping rule.
To deal with such a situation, one has to rely on anytime-valid methods, typically through the use of a test supermartingale (see \cite{HRMS20} and \cite{RRLK20})
or through e-values (see \cite{S20} and \cite{VW20}).
Moreover, e-values are nicer to combine (e.g., using average and product as in \cite{VW20}) especially under arbitrary dependence, in contrast to the complicated methods of merging p-values.

\subsection*{R code}

 An \texttt{R} package \texttt{pmerge} for various merging methods in this paper     is available at \url{https://github.com/YuyuChen-UW/pmerge}.

\subsection*{Acknowledgements}
The authors thank Aaditya Ramdas and Vladimir Vovk for helpful advice on an earlier version of the paper. {The authors thank an Associate Editor and two anonymous reviewers for valuable comments and suggestions.}
Ruodu Wang acknowledges financial support from the Natural Sciences and Engineering Research Council of Canada (RGPIN-2018-03823, RGPAS-2018-522590) and the University of Waterloo  CAE Research Grant  from the Society of Actuaries.

\newpage
\appendix
\renewcommand{\theequation}{\Alph{section}.\arabic{equation}}
\renewcommand{\thetable}{\Alph{section}.\arabic{table}}
\renewcommand{\theexample}{\Alph{section}.\arabic{example}}

\begin{center}
\Large Supplementary Material for \\Trade-off between validity and efficiency of merging p-values under arbitrary dependence
\end{center}
  \section{Proofs of theorems and propositions}\label{section1}
\subsection{Proof of Proposition \ref{p1}}
By definition, we have
\begin{align*}
a_F(\epsilon)=\inf\{q_{\epsilon}(F(U_{1},\dots,U_{K}))\mid U_{1},\dots,U_{K}\in\mathcal{U}\},~\epsilon\in (0,1).
\end{align*}
We shall show \begin{align}\label{a3}a_F(\epsilon)=\inf\{q_1(F(V_{1},\dots,V_{K}))\mid V_{1},\dots,V_{K}\in\mathcal{U}_\epsilon\},~\epsilon\in (0,1),\end{align}
where $\mathcal{U}_\epsilon$ denotes the collection of all uniform random variables distributed on $[0,\epsilon]$.  Denote by $S=F(U_{1},\dots,U_{K})$ and $G_S^{-1}(t)=q_{t}(S),~ t\in (0,1]$. We can find $U_S\in\mathcal{U}$ such that $G_S^{-1}(U_S)=S$ a.s. (e.g., Lemma A.32 of \cite{FS16}). Let
$f_i(t)=\mathbb{P}\left(U_i\leq t|U_S<\epsilon\right),~t\in [0,1]$. Then $f_i(U_i)$ conditionally on $U_S<\epsilon$ is a uniform random variable  on $[0,1]$ and $V_i^\epsilon:=\epsilon f_i(U_i)$ conditionally on $U_S<\epsilon$ is a uniform random variable on $[0,\epsilon]$.
We construct the following two random variables:
\begin{align}
S_1=S\id_{\{U_S<\epsilon\}} +d\id_{\{U_S\geq \epsilon\}},~S_2=F(V_1^\epsilon,\dots,V_n^\epsilon)\id_{\{U_S<\epsilon\}}+d\id_{\{U_S\geq \epsilon\}},
\end{align}
where $d>F(\epsilon,\dots,\epsilon)$. Noting the fact that $\epsilon f_i(t)=\mathbb{P}(U_i\leq t, U_S<\epsilon)\leq t,~t\in [0,1]$ and $F$ is increasing, we have $S_1\geq S_2$. Hence $q_\epsilon(S_1)\geq q_\epsilon(S_2)$. Moreover, direct calculation shows $q_\epsilon(S)=q_\epsilon(S_1)$. Thus $q_\epsilon(S)\geq q_\epsilon(S_2) $. Let $\hat{V}_1, \dots, \hat{V}_n$  be uniform random variables on $[0,\epsilon]$ such that $(\hat{V}_1, \dots, \hat{V}_n)$ has the joint distribution identical to the conditional distribution of $(V_1^\epsilon, \dots, V_n^\epsilon)$ on $U_S<\epsilon$. Hence, for $x<d$,
\begin{align*}
\mathbb{P}(S_2\leq x)&=\mathbb{P}(F(V_1^\epsilon,\dots,V_n^\epsilon)\leq x, U_S<\epsilon)\\
&=\epsilon\mathbb{P}(F(V_1^\epsilon,\dots,V_n^\epsilon)\leq x| U_S<\epsilon)\\
&=\epsilon \mathbb{P}(F(\hat{V}_1,\dots,\hat{V}_n)\leq x).
\end{align*}
This implies $q_\epsilon(S_2)= q_1(F(\hat{V}_1,\dots,\hat{V}_n))$. Thus we have
\begin{align*}
a_F(\epsilon)\geq \inf\{q_1(F(V_{1},\dots,V_{K}))\mid V_{1},\dots,V_{K}\in\mathcal{U}_\epsilon\}.
\end{align*}
We next show ``$\leq$" in (\ref{a3}). Take $V_1,\dots, V_n\in \mathcal{U}_\epsilon$ and $U\in\mathcal{U}$ such that $U$ is independent of $V_1,\dots, V_n$. Let $\hat{U}_i=V_i\id_{\{U<\epsilon\}}+U\id_{\{U\geq\epsilon\}},~ i=1,2,\dots, n$. It is clear that $\hat{U}_i\in\mathcal{U},~i=1,2,\dots,n$ and $F(\hat{U}_1,\dots,\hat{U}_n)=F(V_1,\dots,V_n)\id_{\{U<\epsilon\}}+F(U,\dots,U)\id_{\{U\geq\epsilon \}}$. Noting that $F$ is increasing, we have  $q_1(F(V_1,\dots, V_n))=q_\epsilon(F(\hat{U}_1,\dots,\hat{U}_n))$. This implies \begin{align*}
a_F(\epsilon)\leq \inf\{q_1(F(V_{1},\dots,V_{K}))\mid V_{1},\dots,V_{K}\in\mathcal{U}_\epsilon\}.
\end{align*} Therefore, (\ref{a3}) holds. By (\ref{a3}) and the homogeneity of $F$ we have that for  $\epsilon\in (0,1)$,
\begin{align*}
a_F(\epsilon)&=\inf\{q_1(F(V_{1},\dots,V_{K}))\mid V_{1},\dots,V_{K}\in\mathcal{U}_\epsilon\}\\
&=\inf\{q_1(F(\epsilon U_{1},\dots,\epsilon U_{K}))\mid U_{1},\dots,U_{K}\in\mathcal{U}\}\\
&=\epsilon\inf\{q_1(F(U_{1},\dots,U_{K}))\mid U_{1},\dots,U_{K}\in\mathcal{U}\}.
\end{align*}
This completes the proof.
\qed 

 \subsection{Proof of Proposition \ref{prop:bonf}}
It is well known that the   Bonferroni correction yields $a_{F}(\epsilon)=\epsilon/K$.
Also, since the average of identical objects is itself,
$c_{F}(\epsilon)=\epsilon$ for any averaging method, including the Bonferroni method.
 For iid standard uniform random variables $V_{1},\dots,V_{K}$, we have
  $\mathbb{P}(\min\{ V_1,\dots,V_K\}\leq x)=1-(1-x)^{K}$.  Therefore, $b_{F}(\epsilon)=1-(1-\epsilon)^{1/{K}}$  for $\epsilon\in(0,1)$.
 \qed

 \subsection{Proof of Proposition \ref{p3}}
\begin{enumerate}[(a)]
\item
Suppose $r<0$. We first fix  $K$ and  find the asymptotic of $b_{r}$ as $\epsilon\downarrow 0$ satisfying
\begin{align*}
    \mathbb{P}\left(\sum_{i=1}^{K}P_{i}^{r}\geq K\left(b_{r}(\epsilon)\right)^{r}\right) = \epsilon.
\end{align*}
Observe that the random variables  $P_{i}^{r}$, $i=1,\dots,K$, follow a common Pareto distribution with cdf $\mathbb{P}(P_{i}^{r}\leq x)=1-x^{1/r}$, $x\in(1,\infty)$,  $i=1,\dots,K$.  Note that the tail probability of the sum of iid Pareto random variables is asymptotically the same as that of the maximum of the iid Pareto random variables (e.g., \cite{embrechts2013modelling}, Corollary 1.3.2). Hence
\begin{align*}
    \lim_{\epsilon\downarrow 0}\frac{\mathbb{P}\left(\sum_{i=1}^{K}P_{i}^{r}\geq K\left(b_{r}(\epsilon)\right)^{r}\right)}{\mathbb{P}\left(\max\{P_{1}^{r},\dots,P_{K}^{r}\}> K\left(b_{r}(\epsilon)\right)^{r}\right)}=\lim_{\epsilon\downarrow 0}\frac{\epsilon}{1-\left(1-K^{\frac{1}{r}}b_{r}(\epsilon)\right)^K}=1.
\end{align*}
This implies
$$  b_{r}  (\epsilon)\sim  \frac{1-(1-\epsilon)^{\frac{1}{K}}}{K^{\frac{1}{r}}}\sim K^{-1-1/r} \epsilon, \mbox{~~~~as $\epsilon \downarrow 0$}.$$ The case  $K\to\infty$ follows directly from the generalized central limit theorem (e.g.,
Theorem 1.8.1 of \cite{samoradnitsky2017stable}).
\item
If $r=0$, in  a similar way, we first have,
\begin{align*}
    \mathbb{P}\left(2\sum_{i=1}^{K}\log\frac{1}{P_{i}}\geq 2K\log\frac{1}{b_{r}(\epsilon)}\right) = \epsilon.
\end{align*}
The random variable  $\log\frac{1}{P_{i}}$, $i=1,\dots,K$, follows exponential distribution with parameter 1. Thus $2\sum_{i=1}^{K}\log\frac{1}{P_{i}}$ follows a chi-square distribution with parameter $2K$. We denote
$q_\alpha (\chi^2_\nu)$  the $\alpha$-quantile of the chi-square distribution with  $\nu$ degrees of freedom. Hence
\begin{align*}
    b_{r}(\epsilon)=\exp\left(-\frac{1}{2K}q_{1-\epsilon} \left(\chi^2_{2K}\right)\right).
\end{align*}

\item
If $r>0$, using the result of \cite{W05}, we have for $0\leq x\leq K^{-r}$,
\begin{align*}
\mathbb{P}\left(M_{r,K}(U_1,\dots,U_K)\leq x\right)&=\mathbb{P}\left(\sum_{i=1}^KU_i^r\leq Kx^r\right)\\
&=\lambda\left\{(x_1,\dots,x_K): \sum_{i=1}^{K}x_i^r\leq Kx^r,~ x_1,\dots,x_K\geq 0\right\} \\
&=\frac{(\Gamma(1+1/p))^K}{\Gamma(1+K/p)}K^{K/r}x^K ,
\end{align*}
where $\lambda$ is the Lebesgue measure.
This implies that if $\epsilon \leq \frac{(\Gamma(1+1/p))^K}{\Gamma(1+K/p)}$,
\begin{align}\label{br}
 b_{r}  (\epsilon)= \frac{(\Gamma(1+K/p))^{1/K}\epsilon^{1/K}}{K^{1/r}\Gamma(1+1/p)}.
\end{align}
The asymptotic behaviour of $b_{r}(\epsilon)$
for fixed $\epsilon\in (0,1)$ as $K\to\infty$ can be obtained by the Central Limit Theorem.
Note that the random variables  $P_{i}^{r}$, $i=1,\dots,K$, follow   a common Beta distribution with mean and variance given by, respectively,  $$\mu=(r+1)^{-1},~\mbox{and}~\sigma^{2}=r^2(1+2r)^{-1}(1+r)^{-2}.$$
The Central Limit Theorem gives $(\sum_{i=1}^{K}P_{i}^{r}-K\mu )/\sqrt{K}\sigma\overset{\rm d}{\rightarrow} \mathrm {N}(0,1)$. Hence
$$b_{r}(\epsilon)\sim \left(\frac{\sigma}{\sqrt{K}}\Phi^{-1}(\epsilon)+\mu\right)^{\frac{1}{r}},~\text{as}~K\to\infty,$$
where $\Phi^{-1}$ is the inverse of the   standard normal distribution function. \qed
\end{enumerate}

\subsection{Proof of Proposition \ref{p4}}
By symmetry of the standard Cauchy distribution,
\begin{align*}
a_{F}(\epsilon)&=\mathcal C \left(\inf\left\{q_{\epsilon}\left(\frac{1}{K}\sum_{i=1}^{K} \mathcal C^{-1}(U_{i})\right)\mid U_1,\dots,U_K\in \mathcal U\right\}\right)\\
&=\mathcal C \left(\frac{-1}{K}\sup\left\{q_{1-\epsilon}\left(\sum_{i=1}^{K} \mathcal C^{-1}(U_{i})\right)\mid U_1,\dots,U_K\in \mathcal U\right\}\right).
\end{align*}
Moreover, $\mathcal C^{-1}(U_i),~ i=1,\dots, K,$ follow the standard Cauchy distribution with decreasing density on $[\mathcal C^{-1}(1-\epsilon),\infty]$ for $\epsilon\in (0,1/2)$.
The proposition follows directly from
applying  Corollary 3.7 of  \cite{WPY13}. \qed

 \subsection{Proof of Theorem \ref{Proic-balance}}

 \begin{enumerate}[(i)]
     \item
      IC-balance of $M_{\phi,K}$ for all $K\in\{2,3,\dots\}$ is equivalent to $\frac{1}{K}\sum_{i=1}^{K}\phi(V_{i})\laweq \phi(U)$ for all $K\in\{2,3,\dots\}$, which is further equivalent to the fact that $\phi(U)$ follows a strictly 1-stable distribution.  We know that strictly 1-stable distributions are Cauchy distributions (see, e.g., Theorem 14.15 of \cite{satolevy}).
  This proves the statement of part (i).
  \item
For the Simes function $S_{\alpha,K}=S_K$, $\alpha_{i}=i$ for $i\in\{1,\dots,K\}$ and $b_{F}(x)=c_{F}(x)=x$ for $x\in[0,1]$. Therefore, $S_{\alpha,K}$ is IC-balanced.

Below we show the opposite direction of the statement.
For $n\in\{2,\dots, K\}$, let $V_{(1)},\dots,V_{(n)}$ be the order statistics for $n$ independent standard uniform random variables $V_{1},\dots, V_{n}$.
Let
$(X_1,\dots,X_{n-1})=(V_{(1)}/V_{(n)},\dots,V_{(n-1)}/V_{(n)})$ which is identically distributed as the order statistics for $n-1$ independent standard uniform random variables, independent of $V_{(n)}$. Hence, for $x\in(0,1/\alpha_{n})$,
\begin{align}\label{recursive}
&   \mathbb{P}\left(S_{\alpha,n}(V_{1},\dots, V_{n})> x\right) \notag \\ &=\mathbb{P}\left(V_{(1)}>x\alpha_{1},\dots,V_{(n-1)}>x\alpha_{n-1},V_{(n)}>x\alpha_{n}\right)\nonumber
\\&=\mathbb{P}\left(X_1>x\alpha_{1}/V_{(n)},\dots,X_{n-1}>x\alpha_{n-1}/V_{(n)},V_{(n)}>x\alpha_1\right)
\nonumber\\
   &=\int_{x\alpha_{n}}^{1}\mathbb{P}\left( X_1>x\alpha_{1}/p,\dots,X_{n-1}>x\alpha_{n-1}/p\right)np^{n-1}\d p\nonumber\\
   &=\int_{x\alpha_{n}}^{1}  \mathbb{P}\left(S_{\alpha,n-1}(V_{1},\dots, V_{n-1})> x/p\right)np^{n-1}\d p,
\end{align}
where  for simplicity we use
$S_{\alpha,n-1}$ for $S_{(\alpha_1,\dots,\alpha_{n-1}),n-1}$.
Note that
\begin{align}\label{recursive1}\mathbb{P}\left(S_{\alpha,1}(V_{1})> x\right)=1-\alpha_1x, ~~~x\in (0,1/\alpha_1).
\end{align}
Plugging (\ref{recursive1}) in (\ref{recursive}), we obtain that $\mathbb{P}\left(S_{\alpha,2}(V_{1}, V_{2})> x\right)$ is a polynomial function of $x$ of degree less than or equal to $2$.
 Recursively, using (\ref{recursive}) we are able to show that the function $\mathbb{P}\left(S_{\alpha,n}(V_{1},\dots, V_{n})> x\right)$ for $x\in(0,1/\alpha_{n})$ is a polynomial of $x$ of degree less than or equal to $n$ for $n=2, \dots, K$. Hence, there exist $K$ constants $\beta_0,\dots,\beta_{K-1}$ such that \begin{align*}\mathbb{P}\left(S_{\alpha,K-1}(V_{1},\dots, V_{K-1})> x\right)=\sum_{i=0}^{K-1}\beta_{i}x^i,~~x\in(0,1/\alpha_{K-1}).\end{align*}
Moreover, noting that $S_{\alpha,K}$ is IC-balanced,
 we have
\begin{equation*}
    \int_{x\alpha_{K}}^{1}  \mathbb{P}\left(S_{\alpha,K-1}(V_{1},\dots, V_{K-1})> x/p\right)Kp^{K-1}\d p=\mathbb{P}\left(S_{\alpha,K}(U,\dots, U)> x\right)=1-x\alpha_{K},
\end{equation*}
for $x\in(0,1/\alpha_{K})$.
 Therefore, we have
\begin{align*}
    \int_{x\alpha_{K}}^{1}  \left(\sum_{i=0}^{K-1}\beta_ix^{i}p^{-i}\right)Kp^{K-1}\d p&=1-x\alpha_{K},
  \end{align*}
  which implies that for $x\in(0,1/\alpha_{K})$,
  \begin{align*}
   \sum_{i=0}^{K-1}\frac{K\beta_i}{K-i}x^i-\left(\sum_{i=0}^{K-1}\frac{K\beta_i}{K-i}\alpha_K^{K-i}\right)x^K =1-x\alpha_{K}.
\end{align*}
Solving the above equation, we get $\beta_{0}=1$, $\beta_{1}=-\frac{K-1}{K}\alpha_{K}$ and $\beta_2=\dots=\beta_{K-1}=0$. Consequently,
\begin{align*}
\mathbb{P}\left(S_{\alpha,K-1}(V_{1},\dots, V_{K-1})> x\right)=1-\frac{K-1}{K}\alpha_Kx, ~~~x\in(0,1/\alpha_{K-1}).
\end{align*} Recursively, using (\ref{recursive}) we have \begin{align}\label{recursive2}\mathbb{P}\left(S_{\alpha,n}(V_{1},\dots, V_{n})> x\right)=1-\frac{n}{K}\alpha_Kx, ~~~x\in(0,1/\alpha_{n})
\end{align} for $n=1,\dots, K$,
which gives, using \eqref{recursive1},
\begin{align}\label{ak}\alpha_K=K\alpha_1.
\end{align} Inserting (\ref{recursive2}) into (\ref{recursive}), we obtain, for $x\in (0,1/\alpha_n)$ and $n=2,\dots,K$,
\begin{align*}
1-\frac{n}{K}\alpha_Kx&=\int_{x\alpha_{n}}^{1}\left(1-\frac{n-1}{K}\alpha_Kx p^{-1}\right)  np^{n-1}\d p \\ & =1-\frac{n}{K}\alpha_Kx+\left(\frac{n}{K}\alpha_K\alpha_n^{n-1}-\alpha_n^n\right)x^n.
\end{align*}
Consequently,
\begin{align*}
  \alpha_n=\frac{n}{K}\alpha_K,~~~n=2,\dots,K,
\end{align*}
which together with (\ref{ak})
 implies $\alpha_n=n\alpha_1, ~ k=1,\dots,K.
 $
This gives the desired statement. \qed
 \end{enumerate}
 {In the following example, we shall employ several theorems from \cite{satolevy}. To make our paper more self-contained, we display the useful part of these theorems as below.

  Theorem 8.1 in \cite{satolevy}: $\mu$ is an infinitely divisible distribution in $\mathbb{R}$ if and only if there exist $d\geq 0$, $\gamma\in\mathbb{R}$ and a measure $\nu$ on $\mathbb{R}$ satisfying $\nu(\{0\})=0$ and $\int_{\mathbb{R}}(|x|^2\wedge 1)\nu(\d x)<\infty$, such that the characteristic function of $\mu$ is
  \begin{align}\label{divisible}\hat{\mu}(z)=\exp\left(-\frac{1}{2} d z^2+i\gamma z+\int_{\mathbb{R}}(e^{izx}-1-izx\id_{[-1,1]}(x))\nu(\d x)\right),~z\in\mathbb{R},\end{align}
  where $\id_{[-1,1]}(\cdot)$ is the indicator function and $i^2=-1$.

 Theorem 27.16 in \cite{satolevy}: Suppose $\mu$ satisfies (\ref{divisible}). If $d=0$ and $\nu$ is discrete with total measure infinite, then $\mu$ is a continuous distribution.
 }
 \begin{example}[IC-balanced generalized mean for a finite $K$]
 \label{ex:1}
 We show that  IC-balance of $M_{\phi,K}$ for a finite $K$  does not imply $M_{\phi,K}$ that $\phi$ is the Cauchy quantile function (up to an affine transform).
 For this purpose, we  construct a continuous distribution $\mu$ such that \begin{align}\label{Stable}\frac{1}{K}\sum_{i=1}^{K}X_{i}\laweq X,
\end{align} where $X$ and $X_i, i=1,\dots, K$ are iid random variables with distribution $\mu$, but $\mu$ is not a Cauchy distribution.
Define
\begin{align*}
\hat{\mu}(z)=\exp\left({{\int_{\mathbb{R}}\left(e^{izx}-1-\id_{[-1,1]}(x)\right)\nu(\d x)}}\right),~ z\in\mathbb{R},
\end{align*}
 where $\nu$ is a symmetric measure on $\mathbb{R}\setminus\{0\}$ satisfying
\begin{align*}
\nu(\{K^{n}\})=\nu(\{-K^{n}\})=K^{-n}, ~n\in\mathbb{Z}, ~\text{and}~\nu\left(\mathbb{R}\setminus\left(\{0\}\cup\bigcup_{n\in \mathbb{Z}}\{K^n, -K^n\}\right)\right)=0.
\end{align*}
It follows from Theorem 8.1 of \cite{satolevy} that $\hat{\mu}$ is the characterization function of some infinitely divisible distribution $\mu$.
Also noting that $\nu(\mathbb{R}\setminus\{0\})=\infty$, by Theorem 27.16 of \cite{satolevy} we know that $\mu$ is a continuous distribution.
By Theorem 14.7 of  \cite{satolevy}, $(\hat{\mu}(z))^b=\hat{\mu}(bz), ~z\in\mathbb{R}, b>0$ holds if and only if
$$T_b\nu(B)=b\nu(B),~\text{and}~\int_{1<|x|\leq b}x\nu(\d x)=0,$$
where $T_b\nu(B)=\nu(b^{-1}B)$ for all Borel sets $B\subset \mathbb{R}$.
By symmetry of $\nu$, $\int_{1<|x|\leq b}x\nu(\d x)=0$ holds for any $b>0$. However, $T_b\nu(B)=b\nu(B)$ holds only for $b\in \{K^n, n\in\mathbb{Z}\}$.
Consequently,
$(\hat{\mu}(z))^b=\hat{\mu}(bz), ~z\in\mathbb{R}$ if and only if $b\in \{K^n, n\in\mathbb{Z}\}$.
This implies that $\mu$ is not a Cauchy distribution (strictly 1-stable distribution) but (\ref{Stable}) holds.
\end{example}

 \subsection{Proof of Theorem \ref{theorem2}}
\begin{enumerate}[(i)]
    \item
     Recall that
\begin{align*}
     \mathcal C^{-1}(x)&=\tan\left(-\frac{\pi}{2} + \pi x \right),~~~x\in(0,1);\\
     \mathcal C(y)&=\frac{1}{\pi}\arctan(y)+\frac{1}{2},~~~y\in \mathbb{R}.
\end{align*}
Note that $\mathcal C^{-1}(x) \sim
-1/(\pi x)$ as $x\downarrow 0$ and
$\mathcal C(y) \sim -1/(\pi y)$ as $y\to- \infty$.
For any $\delta_{1},\delta_{2}\in (0,1/K)$, there exists $0<\epsilon<1$ and $m<0$ such that for all $x\in (0,\epsilon)$ and $y\in(-\infty, m)$,
\begin{align}\label{e12}
-\frac{(1+\delta_{1})}{\pi x}\leq \mathcal C^{-1}(x)\leq -\frac{(1-\delta_1)}{\pi x};
\end{align}
\begin{align}\label{e13}
  -\frac{(1-\delta_2)}{\pi y}\leq \mathcal C(y)\leq  -\frac{(1+\delta_2)}{\pi y}.
\end{align}
For $0<c<1$, there exists $0<\epsilon'<\epsilon$  such that
\begin{align}\label{el4}
\sup_{x\in [\epsilon, c]}\left|\tan\left(-\frac{\pi}{2}+\pi x\right)+\frac{1}{\pi x}\right|\leq \frac{\delta_1}{\pi \epsilon'}.
\end{align}
Take
$(p_1,\dots,p_K)$ such that
$p_{(1)}<\epsilon'$ and
$p_{(K)}\leq c<1$.
 Let   $l= \max \{i=1,\dots,K: p_{(i)}< \epsilon \} $. As a consequence of (\ref{e12}), we have
\begin{align*}
-   \sum_{i=1}^{l} \frac{(1+\delta_1)}{\pi p_{(i)}}\leq \sum_{i=1}^{l}\tan\left(-\frac{\pi}{2}+\pi p_{(i)}\right)\leq -\sum_{i=1}^{l}\frac{(1-\delta_1)}{\pi p_{(i)}}.
\end{align*}
For $j>l$,  \eqref{el4} implies
$$\left|\tan\left(-\frac{\pi}{2}+\pi p_{(j)}\right)+\frac{1}{\pi p_{(j)}}\right|\leq \frac{\delta_1}{\pi \epsilon'}\leq \frac{\delta_1}{\pi p_{(1)}}.$$
Therefore,
\begin{align*}
\sum_{i=1}^{K}\tan\left(-\frac{\pi}{2}+\pi p_{i}\right)&\leq  -\sum_{i=1}^{l}\frac{(1-\delta_1)}{\pi p_{(i)}}-\sum_{i=l+1}^K\frac{1}{\pi p_{(i)}}+\frac{(K-l)\delta_1}{\pi p_{(1)}}\\
&\leq - \sum_{i=1}^{K}\frac{(1-K\delta_1)}{\pi p_{(i)}}
 \\& =-\sum_{i=1}^{K}\frac{(1-K\delta_1)}{\pi p_{i}}.
\end{align*}
Similarly, we can show
\begin{align*}
\sum_{i=1}^{K}\tan\left(-\frac{\pi}{2}+\pi p_{i}\right)
\geq \sum_{i=1}^{K}-\frac{(1+K\delta_1)}{\pi p_{i}}.
\end{align*}
Using (\ref{e13}), for any $(p_1,\dots,p_K)$ satisfying $p_{(1)}<\min( \epsilon',\frac{K\delta_1-1}{K\pi m} )$ and $p_{(K)}\leq c<1$,
$$\frac{1-\delta_2}{1+K\delta_1}M_{-1,K}(p_1,\dots,p_K)\leq M_{\mathcal C,K}(p_{1},\dots,p_{K})\leq \frac{1+\delta_2}{1-K\delta_1}M_{-1,K}(p_1,\dots,p_K).$$
 We establish the claim by letting $\delta_1, \delta_2\downarrow 0$, and the above inequalities hold as long as $p_{(1)}$ is sufficiently small.
\item
The statement
    \begin{align*}
    \mathbb{P}\left(M_{\mathcal C,K}(U_{1},\dots,U_{K})<\epsilon\right)\sim\epsilon~~~\mbox{as~}\epsilon \downarrow 0
    \end{align*}
follows directly from Theorem 1 of \cite{LX19} by noting that standard Cauchy distribution is symmetric at $0$. Below we  show $\mathbb{P}\left(M_{-1,K}(U_{1},\dots,U_{K})<\epsilon\right)\sim \epsilon$ as $\epsilon\downarrow 0$, based on similar techniques as in Theorem 1 of \cite{LX19}. Observe that
$$\mathbb{P}\left(M_{-1,K}(U_{1},\dots,U_{K})<\epsilon\right)=\mathbb{P}\left(\frac{1}{K}\sum_{i=1}^{K}U_{i}^{-1}>1/\epsilon\right).$$
Condition (G) means that for any $1\leq i<j\leq K$, $(\Phi^{-1}(U_{i}),\Phi^{-1}(U_{j}))$ is a bivariate  normal random variable with  $\mathrm {cov}(\Phi^{-1}(U_{i}),\Phi^{-1}(U_{j}))=\sigma_{ij}$, where $\Phi$ is the standard normal distribution function and $\Phi^{-1}$ is its inverse. Clearly, $\sigma_{ij}=1$ implies that $U_i=U_j$ a.s. In this case we can combine them in one and the corresponding coefficient becomes $2/K$.
Thus, it suffices to prove the stronger statement
\begin{align}\label{m-1}\mathbb{P}\left(\sum_{i=1}^{K}w_iU_{i}^{-1}>1/\epsilon\right)\sim\epsilon,~\text{as}~\epsilon\downarrow 0,
\end{align}
where $w_i>0,~ i=1,\dots,K$, $\sum_{i=1}^Kw_i=1$ and $\sigma_{ij}<1,~ i,j=1,\dots,K$.
We choose some positive constant $\delta_{\epsilon}$ depending on $\epsilon$,  such that $\delta_{\epsilon}\to 0$ and $\delta_{\epsilon}/\epsilon\to \infty$ as $\epsilon \downarrow0$.  Denote by  $S=\sum_{i=1}^{K}w_iU_{i}^{-1}$, and define the following events: for $i\in\{1,\dots,K\}$,
\begin{align*}
    A_{i,\epsilon}=\left\{U_{i}^{-1}>\frac{1+\delta_{\epsilon}}{w_i\epsilon}\right\},~~
    B_{i,\epsilon}=\left\{U_{i}^{-1}\leq \frac{1+\delta_{\epsilon}}{w_i\epsilon},~ S>1/\epsilon\right\}.
\end{align*}
Let $A_{\epsilon}=\bigcup_{i=1}^{K}A_{i,\epsilon}$ and $B_{\epsilon}=\bigcap_{i=1}^{K}B_{i,\epsilon}$ and thus we have $$\mathbb{P}\left( S >1/\epsilon\right)=\mathbb{P}(A_{\epsilon})+\mathbb{P}(B_{\epsilon}).$$

First we show $\mathbb{P}(B_{\epsilon})=o(\epsilon)$. Note that $S>1/\epsilon$ implies that there exists $i\in \{1,\dots,K\}$ such that $U_i^{-1}>\frac{1}{w_iK\epsilon}$. Hence,
\begin{align*}
\mathbb{P}\left(B_{\epsilon}\right)&\leq \sum_{i=1}^K\mathbb{P}\left(\frac{1}{w_iK\epsilon}<U_{i}^{-1}\leq \frac{1+\delta_{\epsilon}}{w_i\epsilon}, S>1/\epsilon\right)\\
&\leq \sum_{i=1}^K\mathbb{P}\left(\frac{1}{w_iK\epsilon}<U_{i}^{-1}\leq \frac{1-\delta_{\epsilon}}{w_i\epsilon}, S>1/\epsilon\right)+\sum_{i=1}^K\mathbb{P}\left(\frac{1-\delta_{\epsilon}}{w_i\epsilon}<U_{i}^{-1}\leq \frac{1+\delta_{\epsilon}}{w_i\epsilon}\right)\\
&\leq \sum_{i=1}^K\mathbb{P}\left(\frac{1}{w_iK\epsilon}<U_{i}^{-1}\leq \frac{1-\delta_{\epsilon}}{w_i\epsilon}, S>1/\epsilon\right)+\sum_{i=1}^K w_i\epsilon \left(\frac{1}{1-\delta_\epsilon} -\frac{1}{1+\delta_\epsilon}\right)\\
&=: I_1+I_2.
\end{align*}
Noting that $\delta_\epsilon\downarrow 0$ as $\epsilon\downarrow 0$, we have $I_2=o(\epsilon)$. We next focus on $I_1$. Observe
\begin{align*}
I_1&\leq \sum_{i=1}^K\mathbb{P}\left(\frac{1}{w_iK\epsilon}<U_{i}^{-1}\leq \frac{1-\delta_{\epsilon}}{w_i\epsilon}, \sum_{j\neq i}^{K}w_jU_j^{-1}>\delta_\epsilon/\epsilon\right)\\
&\leq \sum_{i=1}^K \sum_{j\neq i}^{K}\mathbb{P}\left(\frac{1}{w_iK\epsilon}<U_{i}^{-1}\leq \frac{1-\delta_{\epsilon}}{w_i\epsilon}, U_j^{-1}>\frac{\delta_\epsilon}{w_jK\epsilon}\right).
\end{align*}
It remains to show for $1\leq i\neq j\leq K$,
\begin{align*}
I_{i,j}:=\mathbb{P}\left(\frac{1}{w_iK\epsilon}<U_{i}^{-1}\leq \frac{1-\delta_{\epsilon}}{w_i\epsilon}, U_j^{-1}>\frac{\delta_\epsilon}{w_jK\epsilon}\right)=o(\epsilon).
\end{align*}
Condition (G) implies that there exist $Z_{i,j}$ and $\delta_{i,j}$ such that
\begin{align}\label{decom}\Phi^{-1}(U_j)=\sigma_{ij}\Phi^{-1}(U_i)+\delta_{ij}Z_{ij},
\end{align}
where $Z_{ij}$ is a standard normal random variable that is independent of $U_i$ and $\sigma_{ij}^2+\delta_{ij}^2=1$.
If $\sigma_{ij}=-1$, we have $U_i=1-U_j$. This implies that  $I_{i,j}=0$ for $\epsilon>0$ sufficiently small. Next, assume $|\sigma_{ij}|<1$, and write
$\gamma_{ij}=\Phi^{-1}\left(w_iK\epsilon\right)$ if $-1<\sigma_{ij}\leq 0$ and $\gamma_{ij}=\Phi^{-1}\left(\frac{w_i\epsilon}{1-\delta_\epsilon}\right)$ if $0<\sigma_{ij}<1$.
We have
\begin{align*}
I_{i,j}&=\mathbb{P}\left(\frac{1}{w_iK\epsilon}<U_{i}^{-1}\leq \frac{1-\delta_{\epsilon}}{w_i\epsilon}, \sigma_{ij}\Phi^{-1}(U_i)+\delta_{ij}Z_{ij}<\Phi^{-1}\left(\frac{w_jK\epsilon}{\delta_\epsilon}\right)\right)\\
&\leq \mathbb{P}\left(\frac{1}{w_iK\epsilon}<U_{i}^{-1}\leq \frac{1-\delta_{\epsilon}}{w_i\epsilon}, \delta_{ij}Z_{ij}<\Phi^{-1}\left(\frac{w_jK\epsilon}{\delta_\epsilon}\right)-\sigma_{ij}\gamma_{ij}\right)\\
&=\mathbb{P}\left(\frac{1}{w_iK\epsilon}<U_{i}^{-1}\leq \frac{1-\delta_{\epsilon}}{w_i\epsilon}\right) \mathbb{P}\left(\delta_{ij}Z_{ij}<\Phi^{-1}\left(\frac{w_jK\epsilon}{\delta_\epsilon}\right)-\sigma_{ij}\gamma_{ij}\right).
\end{align*}
 Note that $\Phi^{-1}(\epsilon)\sim -\sqrt{-2\ln \epsilon},~\text{as} ~\epsilon\downarrow 0$,
which is a slowly varying function.
Taking $\delta_\epsilon=-1/\log \epsilon$, we have $$\Phi^{-1}\left(\frac{w_i\epsilon}{1-\delta_\epsilon}\right)\sim \Phi^{-1}\left(w_iK\epsilon\right)\sim \Phi^{-1}\left(\frac{w_jK\epsilon}{\delta_\epsilon}\right) \mbox{~~as~}\epsilon\downarrow 0.$$
This implies
\begin{align*}
\Phi^{-1}\left(\frac{w_jK\epsilon}{\delta_\epsilon}\right)-\sigma_{ij}\gamma_{ij}\to-\infty, ~\text{as}~\epsilon\downarrow 0.
\end{align*}
Hence $I_{i,j}=o(\epsilon)$. Consequently, $I_1=o(\epsilon)$ and further $\mathbb{P}(B_{\epsilon})=o(\epsilon)$.
Next, we show $\mathbb{P}(A_{\epsilon})\sim\epsilon$. By the Bonferroni inequality, we have,
$$\sum_{i=1}^{K}\mathbb{P}(A_{i,\epsilon})-\sum_{1\leq i<j\leq K}\mathbb{P}(A_{i,\epsilon}\cap A_{j,\epsilon})\leq\mathbb{P}(A_{\epsilon})\leq\sum_{i=1}^{K}\mathbb{P}(A_{i,\epsilon}).$$
Direct calculation gives
\begin{align*}
  \sum_{i=1}^{K}\mathbb{P}(A_{i,\epsilon})=\sum_{k=1}^K\frac{w_i\epsilon}{1+\delta_\epsilon}\sim\epsilon.
\end{align*}
For any $1\leq i<j\leq K$, since the Gaussian copula is tail independent (e.g.,~Example 7.38 of \cite{MFE15}), we have, writing $w=\max\{w_i,w_j\}$,
\begin{align*}
    \mathbb{P}(A_{i,\epsilon}\cap A_{j,\epsilon})&=  \mathbb{P}\left(U_{i}^{-1}>\frac{1+\delta_{\epsilon}}{w_i\epsilon}, U_{j}^{-1}>\frac{1+\delta_{\epsilon}}{w_j\epsilon}\right)\\
    &\le \mathbb{P}\left(U_{i}<\frac{w \epsilon}{1+\delta_{\epsilon}}, U_{j}<\frac{w \epsilon}{1+\delta_{\epsilon}}\right)
 = o(1) \mathbb{P}\left(U_{1 }<\frac{w \epsilon}{1+\delta_{\epsilon}} \right) = o(1) \epsilon.
    \end{align*}
     Hence $\mathbb{P}(A_{i,\epsilon}\cap A_{j,\epsilon})=o(\epsilon)$. This implies  $\mathbb{P}(A_{\epsilon})\sim\epsilon$, and  we establish \eqref{m-1}.

\item

 By Lemma A.1 of \cite{VW19}, we have
 \begin{align*}
 a_\mathcal H(\epsilon)=\epsilon \left(\sup\left\{q_{0}^+\left(\frac{1}{K}\sum_{i=1}^{K}P_{i}^{-1}\right)\mid P_1,\dots,P_K\in \mathcal U \right\}\right)^{-1},~~\epsilon \in (0,1),
 \end{align*}
 where
 $q_{0}^+(X)=\sup\{x\in \mathbb{R}\mid \mathbb{P}(X\leq x)=0\}$.
Note that for any $\delta>0$, there exists $0<\epsilon_\delta<1$ such that for all $x\in (0,\epsilon_\delta)$
\begin{align*}
-\frac{(1+\delta)}{x}<\tan\left(-\frac{\pi}{2}+x\right)< -\frac{(1-\delta)}{x}.
\end{align*}
For  $\delta>0$, letting $0<\epsilon<\epsilon_\delta/\pi$ and using  Theorem 4.6 in \cite{bernard2014risk}, we have
\begin{align*}
   &\inf\left\{ q_{\epsilon}\left(\frac{1}{K}\sum_{i=1}^{K}\mathcal C^{-1}(P_{i})\right)\mid P_1,\dots,P_K\in \mathcal U\right\}\\
   &=\inf\left\{q_{\epsilon}
   \left(\frac{1}{K}\sum_{i=1}^{K}\tan\left(\pi\left(P_{i}-\frac{1}{2}\right)\right)\right)\mid P_1,\dots,P_K\in \mathcal U \right\}\\
    &= \inf\left\{q_{1}\left(\frac{1}{K}\sum_{i=1}^{K}\tan\left(\pi\left(\epsilon P_{i}-\frac{1}{2}\right)\right)\right)\mid P_1,\dots,P_K\in \mathcal U \right\}\\
    &\leq \inf\left\{q_{1}\left(\frac{1}{K}\sum_{i=1}^{K}-\frac{1-\delta}{\epsilon\pi P_i}\right)\mid P_1,\dots,P_K\in \mathcal U \right\}\\
    &=-\frac{1-\delta}{\epsilon\pi}\sup\left\{q_{0}^+\left(\frac{1}{K}\sum_{i=1}^{K}P_{i}^{-1}\right)\mid P_1,\dots,P_K\in \mathcal U \right\}
 = -\frac{1-\delta}{a_{\mathcal H}(\epsilon)\pi}.
\end{align*}
 Similarly, we obtain, for $0<\epsilon<\epsilon_\delta/\pi$,
\begin{align*}
\inf\left\{ q_{\epsilon}\left (\frac{1}{K}\sum_{i=1}^{K}\mathcal C^{-1}(P_{i})\right)\right\}\geq -\frac{1+\delta}{a_{\mathcal H}(\epsilon)\pi}.
\end{align*}
Consequently,
\begin{align*}
 \inf\left\{ q_{\epsilon}\left( \frac{1}{K}\sum_{i=1}^{K}\mathcal C^{-1}(P_{i})\right)\right\} \sim -\frac{1}{a_{\mathcal H}(\epsilon)\pi} ~~~\mbox{as }{\epsilon\downarrow 0}.
\end{align*}
Plugging  the above result in the formula for $a_{\mathcal C}$ in \eqref{Cauchy}, and using
$\mathcal C(y) \sim -1/(\pi y)$ as $y\to- \infty$,
 we have, as $\epsilon\downarrow 0$,
\begin{align*}
    a_{\mathcal C}(\epsilon)&
    = \mathcal C\left(\inf \left\{
    q_{\epsilon}\left( \frac{1}{K}\sum_{i=1}^{K}\mathcal C^{-1}(P_{i})\right)\right\}
    \right)\\&
  \sim -\frac 1 \pi \left(\inf\left\{
    q_{\epsilon}\left( \frac{1}{K}\sum_{i=1}^{K}\mathcal C^{-1}(P_{i})\right)\right\}
    \right)^{-1} \sim
 a_{\mathcal H}(\epsilon).
\end{align*}
This completes the proof.

\item By (i), it suffices to show that for $r\neq -1$
$$\frac{M_{-1,K}(p_{1},\dots,p_{K})}{M_{r,K}(p_{1},\dots,p_{K})}\nrightarrow 1,~ \text{as}\max_{i\in\{1,\dots,K\}}p_{i}\downarrow 0.$$
Take $p_1=p^2$ and  $p_i=x_ip$ with $x_i>0$ and $p>0$ for $i=2,\dots,K$. By homogeneity of $M_r$, for $r\leq -1$,
$$\frac{M_{-1,K}(p_{1},\dots,p_{K})}{M_{r,K}(p_{1},\dots,p_{K})}=\frac{M_{-1,K}(p,x_2,\dots,x_K)}{M_{r,K}(p,x_2,\dots,x_K)}.$$
Hence
$$\lim_{p\downarrow 0}\frac{M_{-1,K}(p_{1},\dots,p_{K})}{M_{r,K}(p_{1},\dots,p_{K})}=K^{1/r+1}\neq 1,~~~r<-1.$$
This proves  the claim of (iv) for $r<-1.$ The case for $r>-1$ can be argued  similarly.    \qed
\end{enumerate}

\subsection{Proof of Theorem \ref{sim_har}}

Take arbitrary $p_{1},\dots,p_{K}\in (0,1]$, and
let $j\in\{1,\dots,K\}$ be such that $\min_{k\in\{1,\dots,K\}}p_{(k)}/k=p_{(j)}/j$.  {Noting that $$\sum_{i=1}^{K}\frac{1}{p_{i}}=\sum_{i=1}^{K}\frac{1}{p_{(i)}},~ \text{and}~ \frac{p_{(j)}}{j}\leq \frac{p_{(i)}}{i},~i=1,\dots,K,$$ we have
 $$\frac{S_{K}(p_{1},\dots, p_{K})}{M_{-1,K}(p_{1},\dots, p_{K})}
    =\frac{1}{j}p_{(j)}\left(\sum_{i=1}^{K}\frac{1}{p_{i}}\right)
    =\sum_{i=1}^{K}\frac{1}{j}p_{(j)}\frac{1}{p_{(i)}}
    \leq \sum_{i=1}^{K}\frac{1}{i}p_{(i)}\frac{1}{p_{(i)}}
    =\sum_{i=1}^{K}\frac{1}{i}
    =\ell_{K}.$$
    }
Moreover, $$\frac{S_{K}(p_{1},\dots, p_{K})}{M_{-1,K}(p_{1},\dots, p_{K})}
    =\frac{1}{j}p_{(j)}\left(\sum_{i=1}^{K}\frac{1}{p_{(i)}}\right)\geq \frac{1}{j}p_{(j)}\left(\sum_{i=1}^{j}\frac{1}{p_{(j)}}+\sum_{i=j+1}^{K}\frac{1}{p_{(i)}}\right)\geq 1.$$
    Therefore, $M_{-1,K} \le S_K \le \ell_K M_{-1,K}$.
    The two special cases of equalities are straightforward to check.
    \qed

  \subsection{Proof of Proposition \ref{prop:logK}}
  \begin{enumerate}[(i)]
     \item
     Recall that $a_{F}(x)=a_{F}x$ for $x\in (0,1)$. By (i) of Proposition \ref{p3}, we have $b_{F}(\delta)\sim \delta$ as $\delta\downarrow 0$. Hence
 $
\lim_{\delta\downarrow 0} {b_{F}(\delta)}/ {a_{F}(\delta)}=  {1}/{a_F}.
$
     By Proposition 6 of  \cite{VW19}, we have $a_{F}\sim {1}/{\log K}$, as $K\rightarrow \infty$. Consequently,
\begin{align*}
\lim_{\delta\downarrow 0}\frac{b_{F}(\delta)}{a_{F}(\delta)}\sim \log K,~\text{as}~K\rightarrow\infty.
 \end{align*}
 Moreover, for the harmonic averaging method, $c_F(\epsilon)=\epsilon$. This implies ${c_{F}(\epsilon)}/{a_{F}(\epsilon)}=  {1}/{a_F}$. We establish the claim  by the fact $a_{F}\sim1/\log K$, as $K\rightarrow \infty$.
     \item
     By Theorem \ref{theorem2}, we have $a_{\mathcal C}(\delta)\sim a_{\mathcal H}(\delta)$ and $b_{\mathcal C}(\delta)\sim b_{\mathcal H}(\delta)$ as $\delta\downarrow 0$, which together with (i) leads to
 $$ \lim_{\delta\downarrow 0}\frac{b_{\mathcal C}(\delta)}{a_{\mathcal C}(\delta)}  \sim \log K,~\text{as}~K\to\infty.$$
The rest of the statement follows by noting that $c_{\mathcal C}(\delta)=b_{\mathcal C}(\delta)$.
     \item
For the Simes  method, recall that $a_{F}(x)=x/\ell_K$ and $b_F(x)=c_F(x)=x$.
     The claim  follows directly from the fact that $\ell_K=\sum_{k=1}^{K}\frac{1}{k}\sim \log K,$ as $K\rightarrow \infty$.
      \qed
 \end{enumerate}

 \section{Additional tables}\label{section2}
 In Tables \ref{t2p} and {\ref{t2p1}} we report numerical results of prices for validity for $\epsilon =0.05$ and $0.0001$, respectively.
 
  \begin{table}[h]
    \caption{$b_{F}(\epsilon)/a_{F}(\epsilon)$ and $c_{F}(\epsilon)/a_{F}(\epsilon)$ for $\epsilon=0.05$ and $K\in \{50,100,200,400\}$}
  \label{t2p}
\centering \small
\begin{tabular}{lcccccccc}
\hline
& \multicolumn{2}{c}{$K=50$} & \multicolumn{2}{c}{$K=100$}& \multicolumn{2}{c}{$K=200$}& \multicolumn{2}{c}{$K=400$}\\
   \hline
                    & $b_{F}/a_{F}$ & $c_{F}/a_{F}$ & $b_{F}/a_{F}$ & $c_{F}/a_{F}$ & $b_{F}/a_{F}$ & $c_{F}/a_{F}$ & $b_{F}/a_{F}$ & $c_{F}/a_{F}$ \\

   {Bonferroni} & 1.025 & 50.000 & 1.026 & 100.000& 1.026 & 200.000& 1.026 & 400.000\\

    {Negative-quartic} & 1.367 & 25.071 & 1.367 & 42.164& 1.368 & 70.911 & 1.368 & 119.257 \\

      {Simes} & 4.499 & 4.499 & 5.187 & 5.187& 5.878 & 5.878 & 6.570 & 6.570 \\

    {Cauchy} & 6.623 & 6.623 & 7.463 & 7.463& 8.274 & 8.274  & 9.055 & 9.055\\

    {Harmonic} & 6.793 & 6.625 & 7.650 & 7.459& 8.485 & 8.273 & 9.306 & 9.072\\

    {Geometric} & 15.679 & 2.718 & 16.874 & 2.718& 17.755 & 2.718  & 18.395 & 2.718\\
   \hline
\end{tabular}
\end{table}

 \begin{table}[h]
 {
    \caption{$b_{F}(\epsilon)/a_{F}(\epsilon)$ and $c_{F}(\epsilon)/a_{F}(\epsilon)$ for $\epsilon=0.0001$ and $K\in \{50,100,200,400\}$}
  \label{t2p1}
\centering \small
\begin{tabular}{lcccccccc}
\hline
& \multicolumn{2}{c}{$K=50$} & \multicolumn{2}{c}{$K=100$}& \multicolumn{2}{c}{$K=200$}& \multicolumn{2}{c}{$K=400$}\\
   \hline
                    & $b_{F}/a_{F}$ & $c_{F}/a_{F}$ & $b_{F}/a_{F}$ & $c_{F}/a_{F}$ & $b_{F}/a_{F}$ & $c_{F}/a_{F}$ & $b_{F}/a_{F}$ & $c_{F}/a_{F}$ \\

   {Bonferroni} & 1.000 & 50.000 & 1.000 & 100.000& 1.000 & 200.000& 1.000 & 400.000\\

    {Negative-quartic} & 1.333 & 25.071 & 1.333 & 42.164& 1.333 & 70.911 & 1.333 & 119.257 \\

      {Simes} & 4.499 & 4.499 & 5.187 & 5.187& 5.878 & 5.878 & 6.570 & 6.570 \\

    {Cauchy} & 6.625 & 6.625 & 7.465 & 7.465& 8.274 & 8.274  & 9.055 & 9.055\\

    {Harmonic} & 6.625 & 6.625 & 7.459 & 7.459& 8.272 & 8.272 & 9.071 & 9.071\\

    {Geometric} & 5416.222 & 2.718 & 6601.414 & 2.718& 7523.231 & 2.718  & 8214.151 & 2.718\\
   \hline
\end{tabular}
}
\end{table}

\end{document}